\documentclass[10pt]{amsart}
\usepackage[T1]{fontenc}
\usepackage{geometry}
\usepackage[latin1] {inputenc}
\usepackage{amsmath}
\usepackage{amsfonts, amssymb, textcomp}
\usepackage[colorlinks=flase, linkcolor=red,urlcolor=green, citecolor=blue]{hyperref}
\usepackage{subeqnarray}

\usepackage{latexsym}
\usepackage{fancyhdr}
\usepackage{longtable}
\usepackage{amsmath, amssymb}
\usepackage{graphicx}

\numberwithin{equation}{section}

\newtheorem{lemma}[subsection]{Lemma}
\newtheorem{theorem}[subsection]{Theorem}
\newtheorem{proposition}[subsection]{Proposition}
\newtheorem{corollary}[subsection]{Corollary}
\newtheorem{definition}[subsection]{Definition}
\newtheorem{example}[subsection]{Example}
\newtheorem{remark}[subsection]{Remark}

\newcommand{\RR}{\mathbb{R}}
\newcommand{\CC}{\mathbb{C}}
\newcommand{\NN}{\mathbb{N}}

\let\on=\operatorname

\newenvironment{demo}[1]{\par\smallskip\noindent{\bf #1.}}{\par\smallskip}

\title[Ultraholomorphic extension theorems in the mixed setting]
{Ultraholomorphic extension theorems in the mixed setting}

\author[J.~Jim\'{e}nez-Garrido, J.~Sanz, and G.~Schindl]{Javier Jim\'{e}nez-Garrido, Javier Sanz and Gerhard Schindl}

\begin{document}
\begin{abstract}
The aim of this work is to generalize the ultraholomorphic extension theorems from V. Thilliez in the weight sequence setting and from the authors in the weight function setting (of Roumieu type) to a mixed framework. Such mixed results have already been known for ultradifferentiable classes and it seems natural that they have ultraholomorphic counterparts. In order to have control on the opening of the sectors in the Riemann surface of the logarithm for which the extension theorems are valid we are introducing new mixed growth indices which are generalizing the known ones for weight sequences and functions. As it turns out, for the validity of mixed extension results the so-called order of quasianalyticity (introduced by the second author for weight sequences) is becoming important.
\end{abstract}

\date{\today}

\maketitle
\section{Introduction}
In the authors' recent works \cite{sectorialextensions} and \cite{sectorialextensions1} we have shown extension theorems in the ultraholomorphic weight function framework, in the first article for spaces of Roumieu type and in the second one also for the Beurling type classes. Such results have already been known before for the weight sequence approach, see \cite{Thilliezdivision}. In \cite{sectorialextensions1} we have transferred Thilliez's ideas to the weight function situation (by using ultradifferentiable Whitney extension results) and in \cite{sectorialextensions} we have used complex methods treated by A. Lastra, S. Malek and the second author~\cite{LastraMalekSanzContinuousRightLaplace, Sanzsummability} in the single weight sequence approach.\vspace{6pt}

In the ultradifferentiable setting also Whitney extension results involving two weight sequences $M$ and $N$ and weight functions $\sigma$ and $\omega$ are known in the literature. In the weight sequence case we refer to \cite{ChaumatChollet94} for the Whitney jet mapping and to \cite{surjectivity} for the Borel mapping, in the weight function case see \cite{BonetMeiseTaylorSurjectivity} for the Borel mapping and \cite{Langenbruch94}, \cite{whitneyextensionmixedweightfunction} and \cite{whitneyextensionmixedweightfunctionII} for the general Whitney jet mapping. In our recent paper \cite{mixedramisurj}, which has served as motivation for this article, by involving a ramification parameter $r\in\NN_{>0}$ we have generalized the mixed setting results from \cite{surjectivity} to $r$-ramification classes introduced in \cite{Schmetsvaldivia00}. We have also generalized the Whitney extension results from \cite{ChaumatChollet94} by using a parameter $r>0$ (see \cite[Theorem 5.10]{mixedramisurj}).

The possibility of an extension in these mixed settings has been characterized in terms of growth properties of weight sequences and functions. We refer also to Remarks \ref{historysequ} and \ref{historyfct} where more (historical) explanations will be given.\vspace{6pt}

From this theoretical point of view it seems natural to ask whether in the ultraholomorphic framework we can also prove extension results in the mixed settings and this question will be treated in this present work. We will consider Roumieu type classes in both the weight sequence and weight function setting. By inspecting the proofs of the main results in \cite{sectorialextensions}, \cite{sectorialextensions1} and \cite{Thilliezdivision} it has turned out that, up to our ability, only the complex methods from \cite{sectorialextensions} admit the possibility to generalize the result to a mixed situation, see Remark \ref{notpossible} below for further details.\vspace{6pt}

The existence of ultraholomorphic extension results is tightly connected to the opening of the sectors where the functions are defined. In the previous results \cite{Thilliezdivision}, \cite{injsurj}, \cite{sectorialextensions} and \cite{sectorialextensions1} growth indices $\gamma(M)$ and $\gamma(\omega)$ have been introduced to measure the maximum size of these sectors, for a detailed study and comparison of these values we refer to \cite{firstindexpaper}. Then a similar notion is required in the mixed setting to obtain satisfactory theorems. Therefore, motivated by the occurring mixed ramified conditions between $M$ and $N$ and their associated weight functions $\omega_M$ and $\omega_N$ appearing in \cite{mixedramisurj} the definition of the mixed growth index for sequences $\gamma(M,N)$ and for weight functions $\gamma(\sigma,\omega)$ has been given, see Section \ref{mixedgrowthsub}.

Under restrictions of the opening of the sectors in terms of these indices, we have stated the main extension result, Theorem \ref{theoExtensionOperatorsMatrixmixed}, for a pair of two given weight functions, using the weight matrix tool described and used in \cite{dissertation} and \cite{compositionpaper}. Then the results are transferred to the weight sequence case thanks to the associated weight functions.\vspace{6pt}

Compared with the previous known extension results for weight functions (in the ultraholomorphic setting) we will also treat ''exotic'' cases here, more precisely: The growth property $\omega(2t)=O(\omega(t))$ as $t\rightarrow+\infty$, denoted by \hyperlink{om1}{$(\omega_1)$} in this article, will not be needed in general anymore in the mixed situation. This property is usually a very basic assumption when working with (Braun-Meise-Taylor) weight functions $\omega$ and it is equivalent to having $\gamma(\omega)>0$ as shown by the authors in \cite{sectorialextensions}. Moreover \hyperlink{om1}{$(\omega_1)$} has also been used to have that the class defined by $\omega$ admits a representation by using the so-called associated weight matrix $\Omega$, see Section \ref{weightmatrixfromfunction} for a summary. Our main extension result Theorem \ref{theoExtensionOperatorsMatrixmixed} is formulated between ultraholomorphic classes defined by weight matrices and we are able to treat such a general situation since in \cite{sectorialextensions} we have worked with weight functions and their associated weight matrices also in a ''nonstandard'' setting, i.e. not assuming \hyperlink{om1}{$(\omega_1)$} necessarily. More detailed explanations will be given in Remark \ref{omega1loss} below. In Appendix \ref{examplesfromlangenbruch} such nonstandard examples will be constructed explicitly and underlining the different situation in our work here.\vspace{6pt}

In the preceding extension results for one sequence, the opening of the sector where the functions are defined is at most $\pi\gamma(M)$. As it will be seen in Section \ref{mixedgrowth}, for any sequences $M$ and $N$ satisfying standard assumptions the mixed index $\gamma(M,N)$ is always belonging to the interval $[\gamma(N),\mu(N)]$, where $\mu(N)$ is denoting the so-called {\itshape order of quasianalyticity} introduced by the second author, see \cite{Sanzflatultraholomorphic} and \cite{injsurj}. We know that even for {\itshape strongly regular sequences} $N$ one can have $\gamma(N)<\mu(N)$ and the gap can become as large as desired, see Remark \ref{gammavsomegaremark}. In these situations, we can provide an extension map for any opening $\pi\gamma$ with $\gamma(N)\le\gamma<\mu(N)$ by limiting the size of the derivatives at the origin in terms of a smaller sequence $M$. Furthermore, this sequence $M$ can be chosen optimal in some sense, thanks to a modified version of the technical construction in \cite[Section 4.1]{whitneyextensionweightmatrix}. Hence we can show that the Borel map will be not surjective necessarily anymore but admitting a controlled loss of regularity, so that $\mu(N)$, usually related to the injectivity of the Borel mapping, does have also a meaning associated with the surjectivity. For weight functions the situation is analogous by introducing the order $\mu(\omega)$ in Section \ref{ordersofquasi}.\vspace{6pt}

The paper is organized as follows: First, in Section \ref{section2} all necessary notation and conditions on weight sequences and functions used in this article will be introduced. In Section \ref{mixedgrowth} we will define and study the new mixed growth indices $\gamma(M,N)$ and $\gamma(\sigma,\omega)$ and investigate also the connection of these values to the orders $\mu(N)$ and $\mu(\omega)$. In Sections \ref{optimalflat} and \ref{sectRightInver1var} we will transfer the results from \cite{sectorialextensions} to the mixed settings and providing only the necessary changes in the proofs, the main results will be Theorem \ref{theoExtensionOperatorsMatrixmixed} for the general mixed weight function case, Corollary \ref{coroExtensionOperatorsWeights} for mixed Braun-Meise-Taylor weight functions having \hyperlink{om1}{$(\omega_1)$} and Theorem \ref{theoExtensionOperatorssequencemixed} for the mixed weight sequence case. In Section \ref{localweightsequence} we will prove mixed extension results fixing only the weight that defines the function space for any sector with opening smaller than $\pi\gamma(\cdot)$, see Theorems \ref{theoExtensionOperatorssequencemixedrem} and \ref{theoExtensionOperatorsMatrixmixedrem}. Finally, in the Appendix \ref{examplesfromlangenbruch}, we are providing some (counter-)examples showing $\gamma(M,N), \gamma(\sigma,\omega)>0$, but such that all nonmixed indices $\gamma(\cdot)$ are vanishing, see Theorem \ref{generalizelangenbruch1}.

\subsection{General notation}
Throughout this paper we will use the following notation: We will write $\NN_{>0}=\{1,2,\dots\}$ and $\NN=\NN_{>0}\cup\{0\}$, moreover we put $\RR_{>0}:=\{x\in\RR: x>0\}$.

\section{Ultradifferentiable classes defined by weight sequences and functions}\label{section2}

\subsection{Weight sequences}\label{section21}

$M=(M_k)_k\in\RR_{>0}^{\NN}$ is called a {\itshape weight sequence}, we introduce also $m=(m_k)_k$ defined by $m_k:=\frac{M_k}{k!}$ and $\mu_k:=\frac{M_k}{M_{k-1}}$, $\mu_0:=1$. Similarly we will use this notation for sequences $N, S, L$ as well. $M$ is called {\itshape normalized} if $1=M_0\le M_1$ holds true and which can always be assumed without loss of generality.

For any given weight sequence $M$ and $r>0$ we will write $M^{1/r}:=((M_p)^{1/r})_{p\in\NN}$.

$(1)$ $M$ is called {\itshape log-convex}, we will write \hypertarget{lc}{$(\text{lc})$}, if
$$\forall\;j\in\NN_{>0}:\;M_j^2\le M_{j-1} M_{j+1},$$
equivalently if $(\mu_p)_p$ is nondecreasing. $M$ is called {\itshape strongly log-convex} if \hyperlink{lc}{$(\text{lc})$} holds for the sequence $m$. If $M$ is log-convex and normalized, then $M$ and the mapping $j\mapsto(M_j)^{1/j}$ are nondecreasing, e.g. see \cite[Lemma 2.0.4]{diploma}. In this case we get $M_k\ge 1$ for all $k\ge 0$ and
\begin{equation*}
\forall\;k\in\NN_{>0}:\;\;\;(M_k)^{1/k}\le\mu_k.
\end{equation*}

$(2)$ $M$ has {\itshape moderate growth}, denoted by \hypertarget{mg}{$(\text{mg})$}, if
$$\exists\;C\ge 1\;\forall\;j,k\in\NN:\;M_{j+k}\le C^{j+k} M_j M_k.$$
We can replace in this condition $M$ by $m$ and by $M^{1/r}$ ($r>0$ arbitrary) by changing the constant $C$.


$(3)$ $M$ is called {\itshape nonquasianalytic,} we write \hypertarget{mnq}{$(\text{nq})$}, if
$$\sum_{p=1}^{\infty}\frac{1}{\mu_p}<\infty.$$

More generally, for arbitrary $r>0$ we call $M$ to be $r$-nonquasianalytic, denoted by \hypertarget{mnqr}{$(\text{nq}_r)$}, if
$$\sum_{p=1}^{\infty}\left(\frac{1}{\mu_p}\right)^{1/r}<\infty,$$
and so $M$ has \hyperlink{mnqr}{$(\text{nq}_r)$} if and only if $M^{1/r}$ has \hyperlink{mnq}{$(\text{nq})$}.

$(4)$ $M$ has \hypertarget{gamma1}{$(\gamma_1)$}, if
$$\sup_{p\in\NN_{>0}}\frac{\mu_p}{p}\sum_{k\ge p}\frac{1}{\mu_k}<\infty.$$
In the literature \hyperlink{gamma1}{$(\gamma_1)$} is also called ``strongly nonquasianalyticity condition''.\vspace{6pt}

Due to technical reasons it is often convenient to assume several properties for $M$ at the same time and hence we define the class\vspace{6pt}

\centerline{$M\in\hypertarget{SRset}{\mathcal{SR}}$, if $M$ is normalized and has \hyperlink{slc}{$(\text{slc})$}, \hyperlink{mg}{$(\text{mg})$} and \hyperlink{gamma1}{$(\gamma_1)$}.}\vspace{6pt}

Using this notation we see that $M\in\hyperlink{SRset}{\mathcal{SR}}$ if and only if $m$ is a {\itshape strongly regular sequence} in the sense of \cite[1.1]{Thilliezdivision} (and this terminology has also been used by several authors so far, e.g. see \cite{Sanzflatultraholomorphic}, \cite{Sanzsummability}). At this point we want to make the reader aware that here we are using the same notation as it has already been used by the authors in \cite{sectorialextensions} and \cite{sectorialextensions1}, whereas in \cite{Thilliezdivision} and also in \cite{firstindexpaper} the sequence $M$ is precisely $m$ in the notation in this work. 

$(5)$ For two weight sequences $M=(M_p)_p$ and $N=(N_p)_p$ we write $M\le N$ if and only if $M_p\le N_p\Leftrightarrow m_p\le n_p$ holds for all $p\in\NN$ (and similarly for the sequence of quotients $\mu$ and $\nu$) and write $M\hypertarget{mpreceq}{\precsim}N$ if
$$\exists\;C_1,C_2\ge 1\;\forall\;p\in\NN:\; M_p\le C_1 C_2^p N_p\Longleftrightarrow\sup_{p\in\NN_{>0}}\left(\frac{M_p}{N_p}\right)^{1/p}<+\infty$$
and call them {\itshape equivalent}, denoted by $M\hypertarget{approx}{\approx}N$, if
$$M\hyperlink{mpreceq}{\precsim}N\;\text{and}\;N\hyperlink{mpreceq}{\precsim}M.$$
In the relations above one can replace $M$ and $N$ simultaneously by $m$ and $n$ because $M\hyperlink{mpreceq}{\precsim}N\Leftrightarrow m\hyperlink{mpreceq}{\precsim}n$.\vspace{6pt}

Some properties for weight sequences are very basic and so we introduce for convenience the following set:
$$\hypertarget{LCset}{\mathcal{LC}}:=\{M\in\RR_{>0}^{\NN}:\;M\;\text{is normalized, log-convex},\;\lim_{k\rightarrow\infty}(M_k)^{1/k}=+\infty\}.$$

It is well-known (e.g. see \cite[Lemma 2.2]{whitneyextensionweightmatrix}) that for any $M\in\hyperlink{LCset}{\mathcal{LC}}$ condition \hyperlink{mg}{$(\on{mg})$} is equivalent to $\sup_{p\in\NN}\frac{\mu_{2p}}{\mu_p}<\infty$ and to $\sup_{p\in\NN_{>0}}\frac{\mu_{p+1}}{(M_p)^{1/p}}<\infty$.

A prominent example are the {\itshape Gevrey sequences} $G^r:=(p!^r)_{p\in\NN}$, $r>0$, which belong to the class $\hyperlink{SRset}{\mathcal{SR}}$ for any $r>1$.

\subsection{Weight functions}\label{section22}
A function $\omega:[0,+\infty)\rightarrow[0,+\infty)$ is called a {\itshape weight function} (in the terminology of \cite[Section 2.2]{sectorialextensions}), if it is continuous, nondecreasing, $\omega(0)=0$ and $\lim_{t\rightarrow+\infty}\omega(t)=+\infty$. If $\omega$ satisfies in addition $\omega(t)=0$ for all $t\in[0,1]$, then we call $\omega$ a {\itshape normalized weight function} and for convenience we will write that $\omega$ has $\hypertarget{om0}{(\omega_0)}$.\vspace{6pt}

Given an arbitrary function $\omega$ we will denote by $\omega^{\iota}$ the function $\omega^{\iota}(t):=\omega(1/t)$ for any $t>0$. Moreover, for $r>0$, we put $\omega^r$ to be the function given by $\omega^r(t):=\omega(t^r)$.\vspace{6pt}

Moreover we consider the following conditions, this list of properties has already been used in ~\cite{dissertation}.

\begin{itemize}
\item[\hypertarget{om1}{$(\omega_1)}$] $\omega(2t)=O(\omega(t))$ as $t\rightarrow+\infty$, i.e. $\exists\;L\ge 1\;\forall\;t\ge 0:\;\;\;\omega(2t)\le L(\omega(t)+1)$.

\item[\hypertarget{om2}{$(\omega_2)$}] $\omega(t)=O(t)$ as $t\rightarrow+\infty$.

\item[\hypertarget{om3}{$(\omega_3)$}] $\log(t)=o(\omega(t))$ as $t\rightarrow+\infty$ ($\Leftrightarrow\lim_{t\rightarrow+\infty}\frac{t}{\varphi_{\omega}(t)}=0$, $\varphi_{\omega}$ being the function defined next).

\item[\hypertarget{om4}{$(\omega_4)$}] $\varphi_{\omega}:t\mapsto\omega(e^t)$ is a convex function on $\RR$.

\item[\hypertarget{om5}{$(\omega_5)$}] $\omega(t)=o(t)$ as $t\rightarrow+\infty$.

\item[\hypertarget{om6}{$(\omega_6)$}] $\exists\;H\ge 1\;\forall\;t\ge 0:\;2\omega(t)\le\omega(H t)+H$.


\item[\hypertarget{omnq}{$(\omega_{\text{nq}})$}] $\int_1^{+\infty}\frac{\omega(t)}{t^2}dt<+\infty.$

\item[\hypertarget{omsnq}{$(\omega_{\text{snq}})$}] $\exists\;C>0\;\forall\;y>0: \int_1^{+\infty}\frac{\omega(y t)}{t^2}dt\le C\omega(y)+C$.
\end{itemize}

An interesting example is $\sigma_s(t):=\max\{0,\log(t)^s\}$, $s>1$, which satisfies all listed properties except \hyperlink{om6}{$(\omega_6)$}. It is well-known that the ultradifferentiable class defined by using the weight $t\mapsto t^{1/s}$ coincides with the ultradifferentiable class given by the weight sequence $G^s=(p!^s)_{p\in\NN}$ of index $s>1$.\vspace{6pt}

For convenience we define the sets
$$\hypertarget{omset0}{\mathcal{W}_0}:=\{\omega:[0,\infty)\rightarrow[0,\infty): \omega\;\text{has}\;\hyperlink{om0}{(\omega_0)},\hyperlink{om3}{(\omega_3)},\hyperlink{om4}{(\omega_4)}\},\hspace{20pt}\hypertarget{omset1}{\mathcal{W}}:=\{\omega\in\mathcal{W}_0: \omega\;\text{has}\;\hyperlink{om1}{(\omega_1)}\}.$$
For any $\omega\in\hyperlink{omset0}{\mathcal{W}_0}$ we define the {\itshape Legendre-Fenchel-Young-conjugate} of $\varphi_{\omega}$ by
\begin{equation*}\label{legendreconjugate}
\varphi^{*}_{\omega}(x):=\sup\{x y-\varphi_{\omega}(y): y\ge 0\},\;\;\;x\ge 0,
\end{equation*}
with the following properties, e.g. see \cite[Remark 1.3, Lemma 1.5]{BraunMeiseTaylor90}: It is convex and nondecreasing, $\varphi^{*}_{\omega}(0)=0$, $\varphi^{**}_{\omega}=\varphi_{\omega}$, $\lim_{x\rightarrow+\infty}\frac{x}{\varphi^{*}_{\omega}(x)}=0$ and finally $x\mapsto\frac{\varphi_{\omega}(x)}{x}$ and $x\mapsto\frac{\varphi^{*}_{\omega}(x)}{x}$ are nondecreasing on $[0,+\infty)$. For any $\omega\in\hyperlink{omset1}{\mathcal{W}}$ we can assume w.l.o.g. that $\omega$ is $\mathcal{C}^1$ (see \cite[Lemma 1.7]{BraunMeiseTaylor90}).\vspace{6pt}

Let $\sigma,\tau$ be weight functions, we write $\sigma\hypertarget{ompreceq}{\preceq}\tau$ if $\tau(t)=O(\sigma(t))\;\text{as}\;t\rightarrow+\infty$
and call them equivalent, denoted by $\sigma\hypertarget{sim}{\sim}\tau$, if
$\sigma\hyperlink{ompreceq}{\preceq}\tau$ and $\tau\hyperlink{ompreceq}{\preceq}\sigma$.

Motivated by the notion of a strong weight function given in \cite{BonetBraunMeiseTaylorWhitneyextension}\vspace{6pt}

\centerline{$\omega$ will be called a strong weight, if $\omega\in\hyperlink{omset0}{\mathcal{W}_0}$ and in addition \hyperlink{omsnq}{$(\omega_{\text{snq}})$} is satisfied.}

\vspace{6pt}
Note that for any weight function property \hyperlink{omnq}{$(\omega_{\text{nq}})$} implies \hyperlink{om5}{$(\omega_5)$} because $\int_t^{+\infty}\frac{\omega(u)}{u^2}du\ge\omega(t)\int_t^{+\infty}\frac{1}{u^2}du=\frac{\omega(t)}{t}$.


If $\omega$ satisfies any of the properties \hyperlink{om1}{$(\omega_1)$}, \hyperlink{om3}{$(\omega_3)$}, \hyperlink{om4}{$(\omega_4)$} or \hyperlink{om6}{$(\omega_6)$}, then the same holds for $\omega^r$, but \hyperlink{om5}{$(\omega_5)$}, \hyperlink{omnq}{$(\omega_{\text{nq}})$} or \hyperlink{omsnq}{$(\omega_{\text{snq}})$} might not be preserved.\vspace{6pt}

Concerning condition \hyperlink{omnq}{$(\omega_{\text{nq}})$} we point out that
\begin{equation}\label{omnqr}
\int_1^{\infty}\frac{\omega^r(u)}{u^2}du=\int_1^{\infty}\frac{\omega(u^r)}{u^2}du=\frac{1}{r}\int_1^{\infty}\frac{\omega(v)}{u^2}\frac{dv}{u^{r-1}}=\frac{1}{r}\int_1^{\infty}\frac{\omega(v)}{v^{1+1/r}}dv,
\end{equation}
hence it makes sense to consider the following generalization \hypertarget{omnqr}{$(\omega_{\text{nq}_r})$} (analogously to \hyperlink{mnqr}{$(\text{nq}_r)$}):
$$\int_1^{\infty}\frac{\omega(u)}{u^{1+1/r}}du<+\infty.$$
Then $\omega^r$ has \hyperlink{omnq}{$(\omega_{\text{nq}})$} if and only if $\omega$ has \hyperlink{omnqr}{$(\omega_{\text{nq}_r})$}.

\subsection{Weight matrices}\label{classesweightmatrices}
For the following definitions see also \cite[Section 4]{compositionpaper}.

Let $\mathcal{I}=\RR_{>0}$ denote the index set (equipped with the natural order), a {\itshape weight matrix} $\mathcal{M}$ associated with $\mathcal{I}$ is a (one parameter) family of weight sequences $\mathcal{M}:=\{M^x\in\RR_{>0}^{\NN}: x\in\mathcal{I}\}$, such that
$$\forall\;x\in\mathcal{I}:\;M^x\;\text{is normalized, nondecreasing},\;M^{x}\le M^{y}\;\text{for}\;x\le y.$$
For convenience we will write \hypertarget{Marb}{$(\mathcal{M})$} for this basic assumption on $\mathcal{M}$. We call a weight matrix $\mathcal{M}$ {\itshape standard log-convex,} denoted by \hypertarget{Msc}{$(\mathcal{M}_{\on{sc}})$}, if $\mathcal{M}$ has \hyperlink{Marb}{$(\mathcal{M})$} and
$$\forall\;x\in\mathcal{I}:\;M^x\in\hyperlink{LCset}{\mathcal{LC}}.$$
Moreover, we put $m^x_p:=\frac{M^x_p}{p!}$ for $p\in\NN$, and $\mu^x_p:=\frac{M^x_p}{M^x_{p-1}}$ for $p\in\NN_{>0}$, $\mu^x_0:=1$.

A matrix is called {\itshape constant} if $M^x\hyperlink{approx}{\approx}M^y$ for all $x,y\in\mathcal{I}$.




Let $\mathcal{M}=\{M^x: x\in\mathcal{I}\}$ and $\mathcal{N}=\{N^x: x\in\mathcal{J}\}$ be \hyperlink{Marb}{$(\mathcal{M})$}. We write $\mathcal{M}\hypertarget{Mroumprecsim}{\{\precsim\}}\mathcal{N}$ if
$$\forall\;x\in\mathcal{I}\;\exists\;y\in\mathcal{J}:\;M^x\hyperlink{precsim}{\precsim}N^y,$$
and call them equivalent, denoted by $\mathcal{M}\hypertarget{Mroumapprox}{\{\approx\}}\mathcal{N}$, if
$\mathcal{M}\hyperlink{Mroumprecsim}{\{\precsim\}}\mathcal{N}$ and $\mathcal{N}\hyperlink{Mroumprecsim}{\{\precsim\}}\mathcal{M}$.

\subsection{Weight matrices obtained by weight functions}\label{weightmatrixfromfunction}
We summarize some facts which are shown in \cite[Section 5]{compositionpaper} and will be needed in this work. All properties listed below will be valid for $\omega\in\hyperlink{omset0}{\mathcal{W}_0}$, except \eqref{newexpabsorb} for which \hyperlink{om1}{$(\omega_1)$} is necessary. 

\begin{itemize}
\item[$(i)$] The idea was that to each $\omega\in\hyperlink{omset0}{\mathcal{W}_0}$ we can associate a \hyperlink{Msc}{$(\mathcal{M}_{\text{sc}})$} weight matrix $\Omega:=\{W^l=(W^l_j)_{j\in\NN}: l>0\}$ by\vspace{6pt}

    \centerline{$W^l_j:=\exp\left(\frac{1}{l}\varphi^{*}_{\omega}(lj)\right)$.}\vspace{6pt}

In general it is not clear that $W^x$ is strongly log-convex, i.e. $w^x$ is log-convex, too.

\item[$(ii)$] $\Omega$ satisfies
    \begin{equation}\label{newmoderategrowth}
    \forall\;l>0\;\forall\;j,k\in\NN:\;\;\;W^l_{j+k}\le W^{2l}_jW^{2l}_k.
    \end{equation}
    In case $\omega$ has moreover \hyperlink{om1}{$(\omega_1)$}, $\Omega$ has also
     \begin{equation}\label{newexpabsorb}
     \forall\;h\ge 1\;\exists\;A\ge 1\;\forall\;l>0\;\exists\;D\ge 1\;\forall\;j\in\NN:\;\;\;h^jW^l_j\le D W^{Al}_j.
     \end{equation}

\item[$(iii)$] Equivalent weight functions yield equivalent associated weight matrices.

\item[$(iv)$] \hyperlink{om5}{$(\omega_5)$} holds if and only if $\lim_{p\rightarrow+\infty}(w^l_p)^{1/p}=+\infty$ for all $l>0$.

\end{itemize}

Moreover we have:
\begin{remark}\label{importantremark}
Let $\omega\in\hyperlink{omset0}{\mathcal{W}_0}$ be given, then $\omega$ satisfies
 \begin{itemize}
 \item[$(a)$] \hyperlink{omnq}{$(\omega_{\on{nq}})$} if and only if some/each $W^l$ satisfies \hyperlink{mnq}{$(\on{nq})$},
 \item[$(b)$] \hyperlink{om6}{$(\omega_6)$} if and only if some/each $W^l$ satisfies \hyperlink{mg}{$(\on{mg})$} if and only if $W^l\hyperlink{approx}{\approx}W^n$ for each $l,n>0$. Consequently \hyperlink{om6}{$(\omega_6)$} is characterizing the situation when $\Omega$ is constant.
 \end{itemize}
\end{remark}

\subsection{Associated weight functions $\omega_M$ and $h_M$}\label{section24}
Let $M\in\RR_{>0}^{\NN}$ (with $M_0=1$), then the {\itshape associated function} $\omega_M: \RR_{\ge 0}\rightarrow\RR\cup\{+\infty\}$ is defined by
\begin{equation*}\label{assofunc}
\omega_M(t):=\sup_{p\in\NN}\log\left(\frac{t^p}{M_p}\right)\;\;\;\text{for}\;t>0,\hspace{30pt}\omega_M(0):=0.
\end{equation*}
For an abstract introduction of the associated function we refer to \cite[Chapitre I]{mandelbrojtbook}, see also \cite[Definition 3.1]{Komatsu73}. If $\liminf_{p\rightarrow\infty}(M_p)^{1/p}>0$, then $\omega_M(t)=0$ for sufficiently small $t$, since $\log\left(\frac{t^p}{M_p}\right)<0\Leftrightarrow t<(M_p)^{1/p}$ holds for all $p\in\NN_{>0}$. Moreover under this assumption $t\mapsto\omega_M(t)$ is a continuous increasing function, which is convex in the variable $\log(t)$ and tends faster to infinity than any $\log(t^p)$, $p\ge 1$, as $t\rightarrow+\infty$. $\lim_{p\rightarrow\infty}(M_p)^{1/p}=+\infty$ implies that $\omega_M(t)<+\infty$ for each finite $t$ and which shall be considered as a basic assumption for defining $\omega_M$.

For all $t,r>0$ we get
\begin{equation}\label{omegaMspower}
\omega_M^r(t)=\omega_M(t^r)=\sup_{p\in\NN}\log\left(\frac{t^{rp}}{M_p}\right)=\sup_{p\in\NN}\log\left(\left(\frac{t^p}{(M_p)^{1/r}}\right)^r\right)=r\omega_{M^{1/r}}(t),
\end{equation}
recalling $M^{1/r}_p=(M_p)^{1/r}$.

We collect some well-known properties for $\omega_M$ (e.g. see \cite[Lem. 2.4]{sectorialextensions} and \cite[Lem. 3.2]{sectorialextensions1}).

\begin{lemma}\label{assofuncproper}
Let $M\in\hyperlink{LCset}{\mathcal{LC}}$ be given, then we get:

\begin{itemize}
\item[$(i)$] $\omega_M$ belongs always to \hyperlink{omset0}{$\mathcal{W}_0$}.

\item[$(ii)$] If $M$ satisfies $\hyperlink{gamma1}{(\gamma_1)}$, then $\omega_M$ has \hyperlink{omsnq}{$(\omega_{\on{snq}})$} (and which implies \hyperlink{om1}{$(\omega_1)$}).

\item[$(iii)$] $M$ has \hyperlink{mg}{$(\on{mg})$} if and only if $\omega_M$ has \hyperlink{om6}{$(\omega_6)$}.




\end{itemize}
\end{lemma}

Let $M\in\RR_{>0}^{\NN}$ (with $M_0=1$) and put
\begin{equation}\label{functionhequ1}
h_M(t):=\inf_{k\in\NN}M_k t^k.
\end{equation}
The functions $h_M$ and $\omega_M$ are related by
\begin{equation}\label{functionhequ2}
h_M(t)=\exp(-\omega_M(1/t))=\exp(-\omega_M^{\iota}(t))\;\;\;\forall\;t>0,
\end{equation}
since $\log(h_M(t))=\inf_{k\in\NN}\log(t^kM_k)=-\sup_{k\in\NN}-\log(t^kM_k)=-\omega_M(1/t)$ (e.g. see also \cite[p. 11]{ChaumatChollet94}).



If $M\in\hyperlink{LCset}{\mathcal{LC}}$, then $M$ has \hyperlink{mg}{$(\text{mg})$} if and only if
\begin{equation}\label{mgsquare}
\exists\;C\ge 1\;\forall\;t>0:\;\;\;h_M(t)\le h_M(Ct)^2,
\end{equation}
e.g. see \cite[Lemma 2.4, Remark 2.5]{whitneyextensionweightmatrix}.

Starting with a weight function we recall the following consequence of \eqref{newmoderategrowth}, see again \cite[Remark 2.5]{whitneyextensionweightmatrix}.

\begin{lemma}\label{functionh2}
Let $\omega\in\hyperlink{omset0}{\mathcal{W}_0}$ be given and $\Omega=\{W^l: l>0\}$ the matrix associated with $\omega$. Then we have
\begin{equation}\label{functionh2equ1}
\exists\;A\ge 1\;\forall\;l>0\;\forall\;s>0:\;\;\;h_{W^l}(s)\le h_{W^{2l}}(As)^2.
\end{equation}
\end{lemma}



\subsection{Classes of ultraholomorphic functions}\label{ultraholomorphicroumieu}
We introduce now the classes under consideration in this paper, see also \cite[Section 2.5]{sectorialextensions} and \cite[Section 2.5]{sectorialextensions1}. For the following definitions, notation and more details we refer to \cite[Section 2]{Sanzflatultraholomorphic}. Let $\mathcal{R}$ be the Riemann surface of the logarithm. We wish to work in general unbounded sectors in $\mathcal{R}$ with vertex at $0$, but all our results will be unchanged under rotation, so we will only consider sectors bisected by direction $0$: For $\gamma>0$ we set $$S_{\gamma}:=\left\{z\in\mathcal{R}: |\arg(z)|<\frac{\gamma\pi}{2}\right\},$$
i.e. the unbounded sector of opening $\gamma\pi$, bisected by direction $0$.

Let $M$ be a weight sequence, $S\subseteq\mathcal{R}$ an (unbounded) sector and $h>0$. We define
$$\mathcal{A}_{M,h}(S):=\{f\in\mathcal{H}(S): \|f\|_{M,h}:=\sup_{z\in S, p\in\NN}\frac{|f^{(p)}(z)|}{h^p M_p}<+\infty\}.$$
$(\mathcal{A}_{M,h}(S),\|\cdot\|_{M,h})$ is a Banach space and we put
\begin{equation*}
\mathcal{A}_{\{M\}}(S):=\bigcup_{h>0}\mathcal{A}_{M,h}(S).
\end{equation*}
$\mathcal{A}_{\{M\}}(S)$ is called the Denjoy-Carleman ultraholomorphic class (of Roumieu type) associated with $M$ in the sector $S$ (it is an $(LB)$ space). Analogously we introduce the space of complex sequences
$$\Lambda_{M,h}:=\{a=(a_p)_p\in\CC^{\NN}: |a|_{M,h}:=\sup_{p\in\NN}\frac{|a_p|}{h^{p} M_{p}}<+\infty\}$$
and put $\Lambda_{\{M\}}:=\bigcup_{h>0}\Lambda_{M,h}$. The (asymptotic) {\itshape Borel map} $\mathcal{B}$ is given by
\begin{equation*}
\mathcal{B}:\mathcal{A}_{\{M\}}(S)\longrightarrow\Lambda_{\{M\}},\hspace{15pt}f\mapsto(f^{(p)}(0))_{p\in\NN},
\end{equation*}
where $f^{(p)}(0):=\lim_{z\in S, z\rightarrow 0}f^{(p)}(z)$.\vspace{6pt}

Similarly as for the ultradifferentiable case, we now define ultraholomorphic classes associated with a normalized weight function $\omega$ satisfying \hyperlink{om3}{$(\omega_3)$}. Given an unbounded sector $S$, and for every $l>0$, we first define
$$\mathcal{A}_{\omega,l}(S):=\{f\in\mathcal{H}(S): \|f\|_{\omega,l}:=\sup_{z\in S, p\in\NN}\frac{|f^{(p)}(z)|}{\exp(\frac{1}{l}\varphi^{*}_{\omega}(lp))}<+\infty\}.$$
$(\mathcal{A}_{\omega,l}(S),\|\cdot\|_{\omega,l})$ is a Banach space and we put
\begin{equation*}
\mathcal{A}_{\{\omega\}}(S):=\bigcup_{l>0}\mathcal{A}_{\omega,l}(S).
\end{equation*}
$\mathcal{A}_{\{\omega\}}(S)$ is called the Denjoy-Carleman ultraholomorphic class (of Roumieu type) associated with $\omega$ in the sector $S$ (it is an $(LB)$ space). Correspondingly, we introduce the space of complex sequences
$$\Lambda_{\omega,l}:=\{a=(a_p)_p\in\CC^{\NN}: |a|_{\omega,l}:=\sup_{p\in\NN}\frac{|a_p|}{\exp(\frac{1}{l}\varphi^{*}_{\omega}(lp))}<+\infty\}$$
and put $\Lambda_{\{\omega\}}:=\bigcup_{l>0}\Lambda_{\omega,l}$. So in this case we get the {\itshape Borel map} $\mathcal{B}:\mathcal{A}_{\{\omega\}}(S)\longrightarrow\Lambda_{\{\omega\}}$.\vspace{6pt}

Finally, we recall the ultradifferentiable function classes of {\itshape Roumieu type} defined by a weight matrix $\mathcal{M}$, introduced in \cite[Section 7]{dissertation} and also in \cite[Section 4.2]{compositionpaper}.

Given a weight matrix $\mathcal{M}=\{M^x\in\RR_{>0}^{\NN}: x\in\RR_{>0}\}$ and a sector $S$ we may define ultraholomorphic classes  $\mathcal{A}_{\{\mathcal{M}\}}(S)$  of {\itshape Roumieu type} as
\begin{equation*}
\mathcal{A}_{\{\mathcal{M}\}}(S):=\bigcup_{x\in\RR_{>0}}\mathcal{A}_{\{M^x\}}(S),
\end{equation*}
and accordingly, $\Lambda_{\{\mathcal{M}\}}:=\bigcup_{x\in\RR_{>0}}\Lambda_{\{M^x\}}$.

Let now $\omega\in\hyperlink{omset1}{\mathcal{W}}$ be given and let $\Omega$ be the associated weight matrix defined in Section \ref{weightmatrixfromfunction} $(i)$, then
\begin{equation}\label{equaEqualitySpacesWeightFunctionMatrix}
\mathcal{A}_{\{\omega\}}(S)=\mathcal{A}_{\{\Omega\}}(S)
\end{equation}
holds as locally convex vector spaces. This equality is an easy consequence of \cite[Lemma 5.9 (5.10)]{compositionpaper} (see \eqref{newexpabsorb}) and the way how the seminorms are defined in these spaces. As one also has $\Lambda_{\{\omega\}}=\Lambda_{\{\Omega\}}$, the Borel map $\mathcal{B}$ makes sense in these last classes,
$\mathcal{B}:\mathcal{A}_{\{\Omega\}}(S)\longrightarrow\Lambda_{\{\Omega\}}$.\vspace{6pt}

In any of the considered ultraholomorphic classes, an element $f$ is said to be {\itshape flat} if $f^{(p)}(0) = 0$ for
every $p\in\NN$, that is, $\mathcal{B}(f)$ is the null sequence.

\section{Mixed growth indices for extension results}\label{mixedgrowth}

\subsection{The indices $\gamma(M,N)$ and $\gamma(\sigma,\omega)$}\label{mixedgrowthsub}
First, for $r>0$ we introduce the following condition which will be denoted by \hypertarget{gammar}{$(\gamma_r)$}, see \cite{Schmetsvaldivia00} for $r\in\NN_{>0}$ and \cite[Lemma 2.2.1]{Thilliezdivision} for $r>0$:
$$\sup_{p\in\NN_{>0}}\frac{(\mu_p)^{1/r}}{p}\sum_{k\ge p}\left(\frac{1}{\mu_k}\right)^{1/r}<+\infty.$$
It is immediate that $M$ has \hyperlink{gammar}{$(\gamma_r)$} if and only if $M^{1/r}$ has \hyperlink{gamma1}{$(\gamma_1)$}.
In \cite[Definition 1.3.5]{Thilliezdivision} the growth index $\gamma(M)$ has been introduced (for strongly regular sequences and using a definition which is not based on property \hyperlink{gammar}{$(\gamma_r)$} directly). In \cite[Thm. 3.11, Cor. 3.12]{firstindexpaper} we have shown for $M\in\hyperlink{LCset}{\mathcal{LC}}$ that
$$\gamma(M)=\sup\{r>0: M\;\text{satisfies}\;\hyperlink{gammar}{(\gamma_r)}\}.$$
One can prove that $M$ has \hyperlink{gammar}{$(\gamma_r)$} if and only if $\gamma(M)>r$ (see \cite[Cor. 3.12]{firstindexpaper}).\vspace{6pt}

Let $M, N\in\RR_{>0}^{\NN}$ with $\mu_p\le C\nu_p$ for some $C\ge 1$ and all $p\in\NN$. For $r>0$ we introduce the following growth property: We write \hypertarget{gammarmix}{$(M,N)_{\gamma_r}$} if
\begin{equation*}
\sup_{p\in\NN_{>0}}\frac{(\mu_p)^{1/r}}{p}\sum_{k\ge p}\left(\frac{1}{\nu_k}\right)^{1/r}<+\infty,
\end{equation*}
and the mixed growth index is defined by
\begin{equation*}
\gamma(M,N):=\sup\{r>0: \hyperlink{gammarmix}{(M,N)_{\gamma_r}}\;\;\text{is satisfied}\}.
\end{equation*}
If none condition \hyperlink{gammarmix}{$(M,N)_{\gamma_r}$} holds true, then we put $\gamma(M,N):=0$. It is evident that $\gamma(M,M)=\gamma(M)$ is valid.

\begin{remark}\label{rprime}
Let $M,N\in\hyperlink{LCset}{\mathcal{LC}}$ with $\mu_p\le C\nu_p$ (note that w.l.o.g. we can take $C=1$, otherwise replace $\nu_p$ by $\widetilde{\nu}_p:=C\nu_p)$.

\begin{itemize}
\item[$(i)$] Given $r>0$, for any $0<r'<r$ we see that \hyperlink{gammarmix}{$(M,N)_{\gamma_r}$} implies \hyperlink{gammarmix}{$(M,N)_{\gamma_{r'}}$}, since we can write
$$\frac{(\mu_p)^{1/r'}}{p}\sum_{k\ge p}\left(\frac{1}{\nu_k}\right)^{1/r'}=\frac{(\mu_p)^{1/r}}{p}\sum_{k\ge p}\left(\frac{1}{\nu_k}\right)^{1/r}\left(\frac{\mu_p}{\nu_k}\right)^{(r-r')/(r'r)},$$
and $\mu_p\le C\nu_p\le C\nu_k$ for all $1\le p\le k$.

\item[$(ii)$] Moreover, in \hyperlink{gammarmix}{$(M,N)_{\gamma_r}$} we can equivalently consider $\frac{(\mu_p)^{1/r}}{p}\sum_{k\ge p+1}\left(\frac{1}{\nu_k}\right)^{1/r}$ since
    $$\frac{(\mu_p)^{1/r}}{p}\sum_{k\ge p+1}\left(\frac{1}{\nu_k}\right)^{1/r}\le\frac{(\mu_p)^{1/r}}{p}\sum_{k\ge p}\left(\frac{1}{\nu_k}\right)^{1/r}\le\frac{C^{1/r}}{p}+\frac{(\mu_p)^{1/r}}{p}\sum_{k\ge p+1}\left(\frac{1}{\nu_k}\right)^{1/r}.$$
\end{itemize}
\end{remark}

In order to see how these definitions have been motivated, we are describing next the appearance of such (non-)mixed relations in the literature.

\begin{remark}\label{historysequ}
Condition \hyperlink{gammar}{$(\gamma_1)$} has appeared as (standard) condition $(M3)$ in \cite{Komatsu73} and in \cite{petzsche} where it has been used to characterize the validity of Borel's theorem in the ultradifferentiable weight sequence setting.

Condition \hyperlink{gammarmix}{$(M,N)_{\gamma_1}$} has appeared in the mixed weight sequence situations in \cite{ChaumatChollet94} (for the Whitney jet map) and in \cite{surjectivity} for the Borel map. More precisely in \cite{surjectivity} it has turned out that the characterizing condition is not \hyperlink{gammarmix}{$(M,N)_{\gamma_1}$} directly, but does coincide with this condition whenever $M$ has \hyperlink{mg}{$(\on{mg})$} (as it has been assumed in \cite{ChaumatChollet94}), see also Remark \ref{remarkthm56} below.

In \cite{Schmetsvaldivia00}, condition \hyperlink{gammar}{$(\gamma_r)$} has appeared (for $r\in\NN_{>0}$) and it has also been used by the authors in \cite{injsurj}. In these works \hyperlink{gammar}{$(\gamma_r)$} played a key-role proving extension theorems for ultraholomorphic classes defined by weight sequences since one is working with auxiliary ultradifferentiable-like function classes first defined in \cite{Schmetsvaldivia00}. In \cite[Lemma 2.2.1]{Thilliezdivision} this condition has been introduced for $r>0$ arbitrary and a connection to the value $\gamma(M)$ has been given.

Finally, condition \hyperlink{gammarmix}{$(M,N)_{\gamma_r}$} has appeared in the recent work by the authors \cite{mixedramisurj} (mainly again for $r\in\NN_{>0}$). There we have generalized the results from \cite{surjectivity} to the auxiliary ultradifferentiable-like function classes, moreover in \cite[Theorem 5.10]{mixedramisurj}, we have given a generalization of the ultradifferentiable Whitney extension results from \cite{ChaumatChollet94} involving a ramification parameter $r>0$.
\end{remark}

Now we turn to the weight function situation. Let $\omega$ be a weight function and $r>0$, we write \hypertarget{gammarfct}{$(\omega_{\gamma_r})$} if
$$\exists\;C>0\;\forall\;t\ge 0:\;\;\;\int_1^{\infty}\frac{\omega(tu)}{u^{1+1/r}}du\le C\omega(t)+C$$
holds true. Using this growth property, by \cite[Lemma 2.10, Theorem 2.11]{firstindexpaper} we have
$$\gamma(\omega)=\sup\{r>0: M\;\text{satisfies}\;\hyperlink{gammarfct}{(\omega_{\gamma_r})}\},$$
with $\gamma(\omega)$ denoting the growth index used and introduced in \cite{sectorialextensions}, \cite{sectorialextensions1} (by considering a different growth property of $\omega$ which is not based on \hyperlink{gammarfct}{$(\omega_{\gamma_r})$}). Note also that $\frac{1}{\gamma(\omega)}$ does coincide with the so-called {\itshape upper Matuszewska index}, see \cite[p. 66]{regularvariation}. For a more detailed study of $\gamma(\omega)$ and its connection to the indices studied in \cite{regularvariation} we refer to Section 2 in the authors' recent work \cite{firstindexpaper}.\vspace{6pt}

Let $\omega,\sigma$ be weight functions with $\sigma\hyperlink{ompreceq}{\preceq}\omega$ (i.e. $\omega(t)=O(\sigma(t))$) and $r>0$, we write \hypertarget{gammarfctmix}{$(\sigma,\omega)_{\gamma_r}$} if
\begin{equation*}
\exists\;C>0\;\forall\;t\ge 0:\;\;\;\int_1^{\infty}\frac{\omega(tu)}{u^{1+1/r}}du\le C\sigma(t)+C,
\end{equation*}
and the mixed growth index is defined by
\begin{equation*}
\gamma(\sigma,\omega):=\sup\{r>0: \hyperlink{gammarfctmix}{(\sigma,\omega)_{\gamma_r}}\;\;\text{is satisfied}\}.
\end{equation*}
If none condition \hyperlink{gammarfctmix}{$(\sigma,\omega)_{\gamma_r}$} holds true, then we put $\gamma(\sigma,\omega):=0$ and it is immediate that $\gamma(\omega,\omega)=\gamma(\omega)$. It is also clear that \hyperlink{gammarfctmix}{$(\sigma,\omega)_{\gamma_r}$} implies \hyperlink{gammarfctmix}{$(\sigma,\omega)_{\gamma_{r'}}$} for all $0<r'<r$.

Similarly as before we are describing the appearance of such (non-)mixed relations in the literature.

\begin{remark}\label{historyfct}
\hyperlink{gammarfct}{$(\omega_{\gamma_1})$}, which is precisely \hyperlink{omsnq}{$(\omega_{\on{snq}})$}, has appeared for $\omega=\omega_M$ in \cite{Komatsu73}, and in \cite{BonetBraunMeiseTaylorWhitneyextension} this condition has been characterized in terms of the validity of the ultradifferentiable Whitney extension theorem in the weight function setting.

The mixed condition \hyperlink{gammarfctmix}{$(\sigma,\omega)_{\gamma_1}$} has been treated in \cite{BonetMeiseTaylorSurjectivity} for the Borel map and in \cite{whitneyextensionmixedweightfunction} and \cite{whitneyextensionmixedweightfunctionII} for the general Whitney jet map (see also \cite{Langenbruch94} for compact {\itshape convex} sets). In these works, condition \hyperlink{gammarfctmix}{$(\sigma,\omega)_{\gamma_1}$} has been identified as the characterizing property.

Finally, in \cite[Theorem 5.10]{mixedramisurj} we have introduced \hyperlink{gammarfctmix}{$(\omega_M,\omega_N)_{\gamma_r}$} in order to prove a generalization of the ultradifferentiable Whitney extension results from \cite{ChaumatChollet94} (again by involving a ramification parameter $r>0$).
\end{remark}

The next observation gives the connection between $\gamma(N)$ and $\gamma(M,N)$, resp. between $\gamma(\omega)$ and $\gamma(\sigma,\omega)$.

\begin{lemma}\label{indicescomparison}
Let $M,N\in\hyperlink{LCset}{\mathcal{LC}}$ be given with $\mu_p\le\nu_p$ and $\omega,\sigma$ be weight functions with $\sigma\hyperlink{ompreceq}{\preceq}\omega$. Then we have
$$\gamma(N)\le\gamma(M,N),\hspace{30pt}\gamma(\omega)\le\gamma(\sigma,\omega).$$
\end{lemma}

\demo{Proof}
First, if $\gamma(N)=0$, $\gamma(\omega)=0$, then the conclusion is clear. If these values are strictly positive, then for any $0<r<\gamma(N),\gamma(\omega)$ we get that \hyperlink{gammar}{$(\gamma_r)$} for $N$ and \hyperlink{gammarfct}{$(\omega_{\gamma_r})$} for $\omega$ hold true and so also \hyperlink{gammarmix}{$(M,N)_{\gamma_r}$} and \hyperlink{gammarfctmix}{$(\sigma,\omega)_{\gamma_r}$} are valid (for any $M$ having $\mu_p\le\nu_p$, $\sigma$ having $\sigma\hyperlink{ompreceq}{\preceq}\omega$).
\qed\enddemo

\begin{remark}
The motivation of defining the mixed growth indices (especially for weight functions) was arising by proving \cite[Theorem 5.10]{mixedramisurj}.

First we wish to mention that in \cite[Commentaires 32]{ChaumatChollet94} it was observed (without giving a proof) that there is a connection between \hyperlink{gammarmix}{$(M,N)_{\gamma_1}$} and \hyperlink{gammarfctmix}{$(\sigma,\omega)_{\gamma_1}$} (under suitable basic assumptions on $M,N$): They have stated that \hyperlink{gammarmix}{$(M,N)_{\gamma_1}$} does imply \hyperlink{gammarfctmix}{$(\omega_M,\omega_N)_{\gamma_1}$} and so generalizing \cite[Prop. 4.4]{Komatsu73} to a mixed setting. In \cite[Lemma 5.7]{whitneyextensionmixedweightfunction} a detailed proof of this implication is given, and in \cite[Lemma 5.8]{mixedramisurj} we have generalized this result by involving a ramification parameter $r>0$.
\end{remark}

We recall the next statement which has been shown in \cite[Lemmas 5.8, 5.9]{mixedramisurj} in order to see how $\gamma(M,N)$ and $\gamma(\omega_M,\omega_N)$ are related. This result is the generalization of \cite[Corollary 4.6 (iii)]{firstindexpaper} to the mixed setting.

\begin{lemma}\label{mglemma1}
Let $M,N\in\hyperlink{LCset}{\mathcal{LC}}$ be given with $\mu_p\le\nu_p$ (and which is equivalent to $\mu_p^r\le\nu_p^r$ for all $r>0$ and implies $M^r\le N^r$). Assume that \hyperlink{gammarmix}{$(M,N)_{\gamma_r}$} holds true for $r>0$. Then the associated weight functions are satisfying
\begin{equation*}\label{mglemma1equ}
\exists\;C>0\;\forall\;t\ge 0:\;\;\;\int_1^{\infty}\frac{\omega_N(tu)}{u^{1+1/r}}du\le C\omega_M(t)+C,
\end{equation*}
i.e. \hyperlink{gammarfctmix}{$(\omega_M,\omega_N)_{\gamma_r}$} is satisfied (recall that $M\le N$ implies $\omega_N(t)\le\omega_M(t)$ for all $t\ge 0$).

Consequently, for sequences $M$ and $N$ as assumed above, we always get $\gamma(M,N)\le\gamma(\omega_M,\omega_N)$.

If $M$ does have in addition \hyperlink{mg}{$(\on{mg})$}, then $\gamma(M,N)=\gamma(\omega_M,\omega_N)$ holds true.
\end{lemma}




\begin{remark}\label{remarkthm56}
For the main results about the surjectivity of the Borel map in ramified ultradifferentiable classes \cite[Thm. 3.2, Thm. 4.2, Thm. 5.5]{mixedramisurj} a weaker but characterizing condition has to be considered. More precisely, in \cite{surjectivity} the following condition has been introduced (denoted by $(\ast)$ there):
$$\exists\;s\in\NN_{>0}:\;\;\sup_{p\in\NN_{>0}}\frac{\lambda_{p,s}^{M,N}}{p}\sum_{k=p}^{\infty}\frac{1}{\nu_k}<+\infty,$$
with $\lambda^{M,N}_{p,s}:=\sup_{0\le j<p}\left(\frac{M_p}{s^p N_j}\right)^{1/(p-j)}$.

In \cite{mixedramisurj} we have generalized this condition to
$$\exists\;s\in\NN_{>0}:\;\;\sup_{p\in\NN_{>0}}\frac{(\lambda_{p,s}^{M,N})^{1/r}}{p}\sum_{k=p}^{\infty}\left(\frac{1}{\nu_k}\right)^{1/r}<+\infty,$$
denoted by \hypertarget{SV}{$(M,N)_{\on{SV}_r}$} (and used for $r\in\NN_{>0}$ in the main results). Hence it seems to be reasonable to define
$$\gamma(M,N)_{\on{SV}}:=\sup\{r\in\RR_{>0}: (M,N)_{\on{SV}_r}\;\;\text{is satisfied}\}.$$
In general we only know that \hyperlink{gammarmix}{$(M,N)_{\gamma_r}$} implies \hyperlink{SV}{$(M,N)_{\on{SV}_r}$}, see \cite[Lemma 2.4]{mixedramisurj}. However, under the standard assumptions of the main results, that is $M,N\in\hyperlink{LCset}{\mathcal{LC}}$ with $\mu_p\le\nu_p$ and such that $M$ does have \hyperlink{mg}{$(\on{mg})$}, Lemma \ref{mglemma1} combined with \cite[Lemma 2.4]{mixedramisurj} yield
$$\gamma(\omega_M,\omega_N)=\gamma(M,N)=\gamma(M,N)_{\on{SV}}.$$
\end{remark}

We summarize several more properties.

\begin{remark}\label{usefulremarks}
\begin{itemize}
\item[$(i)$] By definition it is clear that \hyperlink{gammarmix}{$(M,N)_{\gamma_r}$} if and only if \hyperlink{gammarmix}{$(M^{1/r},N^{1/r})_{\gamma_1}$}, recalling $M_p^{1/r}=(M_p)^{1/r}$, $N_p^{1/r}=(N_p)^{1/r}$.

Similarly \hyperlink{gammarfctmix}{$(\sigma,\omega)_{\gamma_r}$} holds true if and only if \hyperlink{gammarfctmix}{$(\sigma^r,\omega^r)_{\gamma_1}$} for the weights $\omega^r$ given by $\omega^r(t)=\omega(t^r)$ and $\sigma^r(t)=\sigma(t^r)$, because $$\int_1^{\infty}\frac{\omega((tu)^r)}{u^2}du=\frac{1}{r}\int_1^{\infty}\frac{\omega(t^rv)}{u^2}\frac{dv}{u^{r-1}}=\frac{1}{r}\int_1^{\infty}\frac{\omega(t^rv)}{v^{1+1/r}}dv,$$
and then replace $t^r$ by $t$ (see also \eqref{omnqr}).

Summarizing, for all $0<r<\gamma(M,N)$ we get \hyperlink{gammarmix}{$(M^{1/r},N^{1/r})_{\gamma_1}$}, whereas for all $0<r<\gamma(\sigma,\omega)$ we have \hyperlink{gammarfctmix}{$(\sigma^r,\omega^r)_{\gamma_1}$} and so
\begin{equation}\label{growthindexpower}
\forall\;r>0:\;\;\;\gamma(M^r,N^r)=r\gamma(M,N),\hspace{30pt}\gamma(\sigma^{1/r},\omega^{1/r})=r\gamma(\sigma,\omega),
\end{equation}
which generalizes this fact from $\gamma(M)$ and $\gamma(\omega)$.

\item[$(ii)$] \hyperlink{gammarfctmix}{$(\sigma,\omega)_{\gamma_r}$} implies $\int_1^{\infty}\frac{\omega(tu)}{u^{1+1/r}}du\ge\omega(t)\int_1^{\infty}\frac{1}{u^{1+1/r}}du=r\omega(t)$ ($\omega$ is nondecreasing).

\item[$(iii)$] The value $p=1$ in \hyperlink{gammarmix}{$(M,N)_{\gamma_r}$} implies that $N^{1/r}$ satisfies \hyperlink{mnq}{$(\on{nq})$}, i.e. $N$ has \hyperlink{mnqr}{$(\on{nq}_r)$}.

The value $t=1$ in \hyperlink{gammarfctmix}{$(\sigma,\omega)_{\gamma_r}$} yields $\int_1^{\infty}\frac{\omega(u)}{u^{1+1/r}}du\le C_1$ and the calculation in $(i)$ proves that $\omega^r$ satisfies \hyperlink{omnq}{$(\omega_{\on{nq}})$}, or \hyperlink{omnqr}{$(\omega_{\on{nq}_r})$} for $\omega$. Consequently \hyperlink{om5}{$(\omega_5)$} follows for $\omega^r$, too.


\end{itemize}
\end{remark}

\subsection{Orders of quasianalyticity $\mu(N)$ and $\mu(\omega)$}\label{ordersofquasi}

In the ultraholomorphic weight sequence setting another important growth index is known and related to the injectivity of the asymptotic Borel map, the so-called {\itshape order of quasianalyticity}. It has been introduced in \cite[Def. 3.3, Thm. 3.4]{Sanzflatultraholomorphic}, see also \cite{JimenezGarridoSanz}, \cite{injsurj} and \cite{firstindexpaper}. We use the notation from \cite{firstindexpaper} to avoid confusion in the weight function case below and to have a unified notation (coming from \cite[p. 73]{regularvariation}).

For given $N\in\hyperlink{LCset}{\mathcal{LC}}$ we set
\begin{equation}\label{orderofquasi}
\mu(N):=\sup\{r\in\RR_{>0}: \sum_{k\ge 1}\left(\frac{1}{\nu_k}\right)^{1/r}<+\infty\}=\sup\{r\in\RR_{>0}:\;N\;\text{has}\;\hyperlink{mnqr}{(\on{nq}_r)}\}=\frac{1}{\lambda_{(\nu_p)_p}},
\end{equation}
with $\lambda_{(\nu_p)_p}:=\inf\left\{\alpha>0: \sum_{p\ge 1}\left(\frac{1}{\nu_p}\right)^{\alpha}<\infty\right\}$ denoting the so-called {\itshape exponent of convergence} of $N$, see \cite[Prop. 2.13, Def. 3.3, Thm. 3.4]{Sanzflatultraholomorphic} and also \cite[p. 145]{injsurj}.

If none \hyperlink{mnqr}{$(\on{nq}_r)$} holds true, then we put $\mu(N):=0$. If $M\in\hyperlink{LCset}{\mathcal{LC}}$ and $M\hyperlink{approx}{\approx}N$, then $\mu(M)=\mu(N)$ since each \hyperlink{mnqr}{$(\on{nq}_r)$} is stable under \hyperlink{approx}{$\approx$}.\vspace{6pt}

A first immediate consequence is the following:

\begin{lemma}\label{gammasmallerthanq}
Let $M,N\in\hyperlink{LCset}{\mathcal{LC}}$ with $\mu_p\le C\nu_p$, then
$$\gamma(M,N)\le\mu(N).$$
\end{lemma}

\demo{Proof}
Note that for $r>\mu(N)$ property \hyperlink{gammarmix}{$(M,N)_{\gamma_r}$} cannot be valid for any choice $M$ (see $(iii)$ in Remark \ref{usefulremarks}).
\qed\enddemo

\begin{remark}\label{gammavsomegaremark}
Lemmas \ref{indicescomparison} and \ref{gammasmallerthanq} together yield
$$\gamma(N)\le\gamma(M,N)\le\mu(N).$$
Hence the mixed index $\gamma(M,N)$ can become crucial whenever $\gamma(N)<\mu(N)$ does hold true. In general the gap between $\gamma(N)$ and $\mu(N)$ can be as large as desired even if $N\in\hyperlink{SRset}{\mathcal{SR}}$, see \cite[Remark 4.13]{firstindexpaper} and \cite[Section 2.2.5]{dissertationjimenez} for more details.
\end{remark}

According to this observation one can ask now the following question: Is it possible to get extension results for values $\gamma>0$ with $\gamma(N)\le\gamma<\mu(N)$? As we will see in Section \ref{sectRightInver1var}, for values $\gamma<\gamma(M,N)$ we can prove extension results in a mixed setting between $M$ and $N$ but it is still not clear how large the gap between $\gamma(M,N)$ and $\mu(N)$ can be in general.\vspace{6pt}

Given $N\in\hyperlink{LCset}{\mathcal{LC}}$ and $r>0$ we consider

\begin{equation}\label{gamma1descendent}
\exists\;M\in\hyperlink{LCset}{\mathcal{LC}},\;\;\mu_p\le C\nu_p:\;\:\;\sup_{p\in\NN_{>0}}\frac{(\mu_p)^{1/r}}{p}\sum_{k\ge p}\left(\frac{1}{\nu_k}\right)^{1/r}<+\infty,
\end{equation}
and
\begin{equation}\label{growthindexdescendent}
\sup\{r\in\RR_{>0}: \eqref{gamma1descendent}\;\;\text{is satisfied}\}.
\end{equation}
If for none $r>0$ \eqref{gamma1descendent} holds true, then the $\sup$ in \eqref{growthindexdescendent} equals $0$. As commented in Remark \ref{rprime}, given $N\in\hyperlink{LCset}{\mathcal{LC}}$ and $r>0$ with having \eqref{gamma1descendent} for some choice $M\in\hyperlink{LCset}{\mathcal{LC}}$, $\mu_p\le C\nu_p$, then this $M$ is sufficient to guarantee \eqref{gamma1descendent} for all $0<r'<r$ as well.

In \cite[Theorem 3.4]{Sanzflatultraholomorphic} (see also \eqref{orderofquasi}) it has been shown that for any $N\in\hyperlink{LCset}{\mathcal{LC}}$ we have
\begin{equation*}\label{orderofquasi1}
\mu(N)=\liminf_{p\rightarrow\infty}\frac{\log(\mu_p)}{\log(p)}.
\end{equation*}
Thus for given $N\in\hyperlink{LCset}{\mathcal{LC}}$ and $0<r<\mu(N)$ we see that $\nu_p\ge p^r$ for all $p\in\NN$ sufficiently large and $C\nu_p\ge p^r$ for all $p\in\NN$ by choosing $C$ large enough. Consequently, in \eqref{gamma1descendent} the choice $M\equiv G^r$, i.e. the Gevrey sequence with index $r>0$, does always make sense and the next result is becoming immediate:

\begin{proposition}\label{gammamixedvsomegam}
Let $N\in\hyperlink{LCset}{\mathcal{LC}}$ be given, then
$$\sup\{r\in\RR_{>0}: \eqref{gamma1descendent}\;\;\text{is satisfied}\}=\mu(N).$$
\end{proposition}

\demo{Proof}
If $\mu(N)=0$, then we have obviously equality. So let now $\mu(N)>0$.

First, if $0<r<\sup\{r\in\RR_{>0}: \eqref{gamma1descendent}\;\;\text{is satisfied}\}$, then the choice $p=1$ in \eqref{gamma1descendent} immediately implies \hyperlink{mnqr}{$(\on{nq}_r)$} for $N$, hence $r\le\mu(N)$ and so the first half is shown.

Conversely, let $r<\mu(N)$ be given, then \eqref{gamma1descendent} is satisfied for $\mu_p=p^r$ and which can be taken as seen above. Hence $\mu(N)\le\sup\{r\in\RR_{>0}: \eqref{gamma1descendent}\;\;\text{is satisfied}\}$ is also shown and we are done.
\qed\enddemo

A disadvantage of taking directly $M\equiv G^r$ is that it is not clear that this precise choice is optimal in the sense that it is the largest sequence $\mu_p\le C\nu_p$ admitting \eqref{gamma1descendent}. To obtain this optimal sequence we recall the following construction: In \cite[Section 4.1]{whitneyextensionweightmatrix}, and which is based on an idea arising in the proof of \cite[Proposition 1.1]{petzsche}, it has been shown that to each $N\in\hyperlink{LCset}{\mathcal{LC}}$ satisfying \hyperlink{mnq}{$(\on{nq})$} we can associate a sequence $S^N$ with good regularity properties and which has been denoted by {\itshape descendant}.\vspace{6pt}

For the reader's convenience we recall now the construction in the following observation and are involving a ramification parameter $r>0$ as well.

\begin{remark}\label{descendant}
Let $N\in\hyperlink{LCset}{\mathcal{LC}}$ be given and satisfying \hyperlink{mnqr}{$(\on{nq}_r)$}, $r>0$. Then there does exist a sequence $S^{N,r}$, the so-called {\itshape descendant of} $N^{1/r}$ defined by its quotients
$$\sigma^{N,r}_k:=\frac{\tau^r_1 k}{\tau^r_k},\hspace{20pt}\sigma^{N,r}_0:=1,$$
with
$$\tau^r_k:=\frac{k}{(\nu_k)^{1/r}}+\sum_{j\ge k}\left(\frac{1}{\nu_j}\right)^{1/r},\;\;\;k\ge 1.$$
So $S^{N,r}$ is depending on both given $N$ and $r>0$ and satisfies the following properties (see \cite[Lemma 4.2]{whitneyextensionweightmatrix}):

\begin{itemize}
\item[$(i)$] $\sigma^{N,r}_k\ge 1$ for all $k\in\NN$, $s^{N,r}:=(S^{N,r}_k/k!)_{k\in\NN}\in\hyperlink{LCset}{\mathcal{LC}}$ (so $S^{N,r}$ is strongly log-convex),

\item[$(ii)$] there exists $C>0$ such that $\sigma^{N,r}_k\le C(\nu_k)^{1/r}$ for all $k\in\NN$,

\item[$(iii)$] \hyperlink{gammarmix}{$(S^{N,r},N^{1/r})_{\gamma_1}$}, then $\gamma(S^{N,r},N^{1/r})\ge 1$,

\item[$(iv)$] if $N$ enjoys \hyperlink{mg}{$(\on{mg})$} (equivalently $N^{1/r}$ does so for some/each $r>0$), then $S^{N,r}$ does have \hyperlink{mg}{$(\on{mg})$} too ($r>0$ arbitrary).

\item[$(v)$] $S^{N,r}$ is optimal/maximal in the following sense: If $M\in\hyperlink{LCset}{\mathcal{LC}}$ is given with $\mu_k\le C(\nu_k)^{1/r}$ and \hyperlink{gammarmix}{$(M,N^{1/r})_{\gamma_1}$}, then $\mu_k\le D\sigma^{N,r}_k$ follows.

We also have that $C^{-1}\sigma^{N,r}_k\le(\nu_k)^{1/r}\le C\sigma_k^{N,r}$ if and only if $N^{1/r}$ does have \hyperlink{gamma1}{$(\gamma_1)$} resp. if and only if  \hyperlink{gammarmix}{$(N,N)_{\gamma_r}$}.
\end{itemize}

We have that $s^{N,r}\le Cs^{N,r'}\Leftrightarrow S^{N,r}\le CS^{N,r'}$ for all $0<r'\le r$ (since $r\mapsto\tau^r_k$ is increasing for all $k\in\NN$ fixed).
\end{remark}

Unfortunately this construction is in general not well-behaved under applying ramification, i.e. $\sigma^{N,r}\neq(\sigma^{N,1})^{1/r}$.

Now we put
\begin{equation}\label{sequenceL}
L^{N,r}:=(S^{N,r})^r.
\end{equation}
Hence $L^{N,r}\in\hyperlink{LCset}{\mathcal{LC}}$ and moreover

\begin{itemize}

\item[$(a)$] \hyperlink{gammarmix}{$((L^{N,r})^{1/r},N^{1/r})_{\gamma_1}$}, equivalently \hyperlink{gammarmix}{$(L^{N,r},N)_{\gamma_r}$} holds true and so $\gamma(L^{N,r},N)\ge r$,

\item[$(b)$] $\lambda^{N,r}_k=(\sigma^{N,r}_k)^r\le C\nu_k$ for all $k\in\NN$ and so \hyperlink{gammarmix}{$(L^{N,r},N)_{\gamma_{r'}}$} for all $0<r'<r$ as well (see $(i)$ in Remark \ref{rprime}),

\item[$(c)$] if $M\in\hyperlink{LCset}{\mathcal{LC}}$ is given with $\mu_k\le C\nu_k$ and \hyperlink{gammarmix}{$(M,N)_{\gamma_r}$}, i.e. \hyperlink{gammarmix}{$(M^{1/r},N^{1/r})_{\gamma_1}$}, then $\mu_k\le D^r(\sigma^{N,r}_k)^r=D^r\lambda^{N,r}_k$ for all $k\in\NN$ and consequently $L^{N,r}$ is maximal (up to a constant) among all sequences satisfying $\mu_k\le C\nu_k$ and \hyperlink{gammarmix}{$(M,N)_{\gamma_r}$},

\item[$(d)$] $S^{N,r}$ has \hyperlink{mg}{$(\on{mg})$} if and only if $L^{N,r}$ does so.

\end{itemize}

According to the maximality $(c)$ mentioned above we point out we have $\lambda^{N,r}_k=(\sigma^{N,r}_k)^r\ge k^r$ for all $k\in\NN$ and so $L^{N,r}\ge G^r$. Moreover, one can show that $\lim_{k\rightarrow\infty}\frac{\lambda^{N,r}_k}{k^r}=\infty$ and which implies $\left(\frac{G^r_k}{L^{N,r}_k}\right)^{1/k}\rightarrow 0$ as $k\rightarrow+\infty$ and so $L^{N,r}$ is strictly larger than $G^r$.\vspace{6pt}



As mentioned above, condition \hyperlink{mg}{$(\on{mg})$} for $N$ does always imply this property for the descendant $S^{N,1}=:S$. However we can obtain a precise characterization of this growth behavior.

\begin{lemma}\label{descendantmg}
Let $N\in\hyperlink{LCset}{\mathcal{LC}}$ having \hyperlink{mnq}{$(\on{nq})$} be given and $S:=S^{N,1}$ shall denote its descendant. Then $S$ does have \hyperlink{mg}{$(\on{mg})$} if and only if
\begin{equation}\label{descendantmgequ}
\exists\;C\ge 1\;\forall\;k\in\NN:\;\;\;\frac{\nu_{2k}}{\nu_k}\le C+C\frac{\nu_{2k}}{2k}\sum_{j\ge 2k}\frac{1}{\nu_j}.
\end{equation}
\end{lemma}

\demo{Proof}
Since $S\in\hyperlink{LCset}{\mathcal{LC}}$, this sequence has \hyperlink{mg}{$(\on{mg})$} if and only if $\sup_{k\in\NN}\frac{\sigma_{2k}}{\sigma_k}<\infty$, e.g. see \cite[Lemma 2.2]{whitneyextensionweightmatrix}. By the definitions given in Remark \ref{descendant} we get that
$\sigma_{2k}\le D\sigma_k\Leftrightarrow\tau_k\le\frac{D}{2}\tau_{2k}$ (with $\tau_k\equiv\tau^1_k$) and which is equivalent to $\frac{k}{\nu_k}+\sum_{j\ge k}\frac{1}{\nu_j}\le D\frac{k}{\nu_{2k}}+\frac{D}{2}\sum_{j\ge 2k}\frac{1}{\nu_j}$ and so to having $\frac{k}{\nu_k}+\sum_{j=k}^{2k-1}\frac{1}{\nu_j}\le D\frac{k}{\nu_{2k}}+(\frac{D}{2}-1)\sum_{j\ge 2k}\frac{1}{\nu_j}$ and so finally \hyperlink{mg}{$(\on{mg})$} is equivalent to
\begin{equation}\label{descendantmgequ1}
\exists\;D\ge 1\;\forall\;k\in\NN:\;\;\;\frac{\nu_{2k}}{\nu_k}+\frac{\nu_{2k}}{k}\sum_{j=k}^{2k-1}\frac{1}{\nu_j}\le D+(\frac{D}{2}-1)\frac{\nu_{2k}}{k}\sum_{j\ge 2k}\frac{1}{\nu_j}.
\end{equation}
The sum on the left-hand side above is estimated by below by $\frac{k}{\nu_{2k}}$ and by above by $\frac{k}{\nu_k}$ (since $(\nu_k)_k$ is increasing). Hence \eqref{descendantmgequ1} implies $\frac{\nu_{2k}}{\nu_k}+1\le D+(D-2)\frac{\nu_{2k}}{2k}\sum_{j\ge 2k}\frac{1}{\nu_j}$ and so \eqref{descendantmgequ} with $C=D-1$ follows, whereas \eqref{descendantmgequ} implies $\frac{\nu_{2k}}{\nu_k}+\frac{\nu_{2k}}{k}\sum_{j=k}^{2k-1}\frac{1}{\nu_j}\le\frac{\nu_{2k}}{\nu_k}+\frac{\nu_{2k}}{\nu_k}\le 2C+C\frac{\nu_{2k}}{k}\sum_{j\ge 2k}\frac{1}{\nu_j}$ and \eqref{descendantmgequ1} follows with $D:=4C$.
\qed\enddemo

Concerning the characterizing condition \eqref{descendantmgequ} we point out that (as expected) it is obviously satisfied if $N$ has \hyperlink{mg}{$(\on{mg})$}, i.e. $\sup_{k\in\NN}\frac{\nu_{2k}}{\nu_k}<\infty$, and moreover:

\begin{itemize}
\item[$(i)$] The right-hand side in \eqref{descendantmgequ} is bounded if and only if \hyperlink{gamma1}{$(\gamma_1)$} is valid for $N$ since for all $k\ge 1$ we get $\frac{\nu_{2k}}{2k}\sum_{j\ge 2k}\frac{1}{\nu_j}\ge\frac{1}{2}\frac{\nu_{2k-1}}{2k-1}\sum_{j\ge 2k-1}\frac{1}{\nu_j}-\frac{1}{2}\frac{\nu_{2k-1}}{2k-1}\frac{1}{\nu_{2k-1}}$ and so $\frac{\nu_{2k-1}}{2k-1}\sum_{j\ge 2k-1}\frac{1}{\nu_j}\le 2\frac{\nu_{2k}}{2k}\sum_{j\ge 2k}\frac{1}{\nu_j}+\frac{1}{2k-1}$.

\item[$(ii)$] In this case $S$ is equivalent to $N$ and so $S$ has \hyperlink{mg}{$(\on{mg})$} if and only if $N$ has this property.

\item[$(iii)$] Instead of having \eqref{descendantmgequ} one can study the more ''compact and easy to handle'' requirement
\begin{equation}\label{descendantmgequalterative}
\liminf_{k\in\NN_{>0}}\frac{\nu_{k}}{k}\sum_{j\ge 2k}\frac{1}{\nu_j}>0.
\end{equation}
Since \eqref{descendantmgequ} is equivalent to having $\frac{2}{C}\le 2\frac{\nu_k}{\nu_{2k}}+\frac{\nu_k}{k}\sum_{j\ge 2k}\frac{1}{\nu_j}$ and \eqref{descendantmgequalterative} to $\frac{\nu_{k}}{k}\sum_{j\ge 2k}\frac{1}{\nu_j}\ge\varepsilon$ for some $\varepsilon>0$ and all $k\in\NN$ we immediately get that \eqref{descendantmgequalterative} implies \eqref{descendantmgequ}.

However, the converse is not clear in general: If $N$ has in addition $\lim_{k\rightarrow\infty}\frac{\nu_{2k}}{\nu_k}=\infty$, i.e. $(\beta^0_2)$ in \cite{petzsche}, then \eqref{descendantmgequ} still can be valid. But, as shown in \cite[Proposition 1.1.(b), $(\beta^0_2)\Rightarrow(\gamma_2)$]{petzsche}, this would yield $\lim_{p\rightarrow\infty}\frac{\nu_{k}}{k}\sum_{j\ge 2k}\frac{1}{\nu_j}=0$ and so making \eqref{descendantmgequalterative} impossible.
\end{itemize}

We are now studying an example which shows that the characterization obtained in Lemma \ref{descendantmg} can hold true even if $N$ does not have \hyperlink{mg}{$(\on{mg})$}.

\begin{example}
Put $\nu_0:=1$ and
$$\nu_k:=2^p((p+1)!-p!)=2^pp!p\hspace{10pt}\text{for}\hspace{10pt}p!\le k<(p+1)!,\;\;\;p\in\NN_{\ge 1}.$$
By definition we have $N\in\hyperlink{LCset}{\mathcal{LC}}$ and $\sum_{j\ge 1}\frac{1}{\nu_j}=\sum_{l=1}^{\infty}\sum_{j=l!}^{(l+1)!-1}\frac{1}{\nu_j}=\sum_{l=1}^{\infty}\frac{(l+1)!-l!}{2^ll!l}=\sum_{l=1}^{\infty}\frac{1}{2^l}=1$, which shows \hyperlink{mnq}{$(\on{nq})$}.

$N$ does not have \hyperlink{mg}{$(\on{mg})$} because for any $k\in\NN_{\ge 1}$ with $k<(p+1)!\le 2k<(p+2)!$, $p\ge 1$, we get $\frac{\nu_{2k}}{\nu_k}=2\frac{(p+1)^2}{p}$ and which tends to infinity as $k\rightarrow\infty$.

Finally let us show that \eqref{descendantmgequalterative} holds (and so \eqref{descendantmgequ}). Let now $k\in\NN$ be given with $p!\le 2k<(p+1)!$, $p\ge 2$. We split the sum
$$\frac{\nu_{k}}{k}\sum_{j\ge 2k}\frac{1}{\nu_j}=\frac{\nu_k}{k}\frac{(p+1)!-2k}{2^pp!p}+\frac{\nu_k}{k}\sum_{l\ge p+1}\frac{(l+1)!-l!}{2^ll!l}=\frac{\nu_k}{k}\frac{(p+1)!-2k}{2^pp!p}+\frac{\nu_k}{k}\frac{1}{2^p},$$
and remark that both summands are nonnegative for all $k\in\NN$ under consideration.

We study the second summand and distinguish between two cases. If $k<p!$, then we have $k\ge p!/2\ge (p-1)!$ and so
$$\frac{\nu_k}{k}\frac{1}{2^p}=\frac{2^{p-1}(p-1)!(p-1)}{k}\frac{1}{2^p}\ge\frac{1}{2}\frac{(p-1)!(p-1)}{p!}=\frac{1}{2}\frac{p-1}{p}\ge\frac{1}{4}.$$
If $p!\le k$, then we can estimate by
$$\frac{\nu_k}{k}\frac{1}{2^p}=\frac{2^pp!p}{k}\frac{1}{2^p}\ge\frac{p!p}{(p+1)!}=\frac{p}{p+1}\ge\frac{1}{2}.$$
Thus the descendant $S$ does have \hyperlink{mg}{$(\on{mg})$}. But since $N$ violates this property, $N$ cannot be equivalent to $S$ and so $N$ does not satisfy \hyperlink{gamma1}{$(\gamma_1)$}.
\end{example}



Similarly, there does exist also an inverse construction concerning the descendant, called the {\itshape predecessor}. However, this does not provide any new insight, see Remark \ref{Javipredecessordestroyer}.\vspace{6pt}

Analogously as in the weight sequence situation we can treat the weight function case as well. First, inspired by the weight sequence approach in \eqref{orderofquasi} and condition \hyperlink{omnqr}{$(\omega_{\text{nq}_r})$}, we define now for any given weight function $\omega$ (as in Section \ref{section22}) the {\itshape order of quasianalyticity} by
\begin{equation}\label{orderofquasifct}
\mu(\omega):=\sup\{r\in\RR_{>0}: \int_1^{\infty}\frac{\omega(u)}{u^{1+1/r}}du<+\infty\}=\sup\{r\in\RR_{>0}:\;\omega\;\text{has}\;\hyperlink{omnqr}{(\omega_{\text{nq}_r})}\}.
\end{equation}

If none condition \hyperlink{omnqr}{$(\omega_{\text{nq}_r})$} holds true, then we put $\mu(\omega):=0$. If $\sigma$ is another weight function with $\omega\hyperlink{sim}{\sim}\sigma\Leftrightarrow\omega^r\hyperlink{sim}{\sim}\sigma^r$, then $\mu(\omega)=\mu(\sigma)$ follows since each condition \hyperlink{omnqr}{$(\omega_{\text{nq}_r})$} is stable under \hyperlink{sim}{$\sim$}.

This growth index likely will have an interpretation for the quasianalyticity of classes of ultraholomorphic functions defined in terms of weight functions, more precisely in order to prove analogous results to \cite{injsurj}, see also \cite{Sanzflatultraholomorphic} and \cite{dissertationjimenez}.\vspace{6pt}

Analogously to \eqref{gamma1descendent} for a given weight function $\omega$ and $r>0$ we consider now
\begin{equation}\label{gamma1omegadescendent}
\exists\;\sigma,\;\sigma\hyperlink{ompreceq}{\preceq}\omega,\;\exists\;C>0\;\forall\;t\ge 0:\;\:\;\int_1^{\infty}\frac{\omega(tu)}{u^{1+1/r}}du\le C\sigma(t)+C,
\end{equation}
and let
\begin{equation}\label{growthindexomegadescendent}
\sup\{r\in\RR_{>0}: \eqref{gamma1omegadescendent}\;\;\text{is satisfied}\}.
\end{equation}
If none \eqref{gamma1omegadescendent} holds true, then the $\sup$ in \eqref{growthindexomegadescendent} equals $0$.

Let $\omega$ and $r>0$ be given and assume that \eqref{gamma1omegadescendent} holds with some $\sigma$, then the same choice is sufficient to have \eqref{gamma1omegadescendent} for all $0<r'<r$ as well.

Lemma \ref{gammasmallerthanq} turns now into:

\begin{lemma}\label{gammasmallerthanqom}
Let $\sigma$ and $\omega$ be weight functions (in the sense of Section \ref{section22}) with $\sigma\hyperlink{ompreceq}{\preceq}\omega$. Then we get
$$\gamma(\sigma,\omega)\le\mu(\omega).$$
\end{lemma}

\demo{Proof}
For any $r>\mu(\omega)$ we see that \hyperlink{omnqr}{$(\omega_{\text{nq}_r})$} is violated and so \hyperlink{gammarfctmix}{$(\sigma,\omega)_{\gamma_r}$} cannot be valid for any choice $\sigma$ (see $(iii)$ in Remark \ref{usefulremarks}).
\qed\enddemo

\begin{remark}
Lemmas \ref{indicescomparison} and \ref{gammasmallerthanqom} together yield
$$\gamma(\omega)\le\gamma(\sigma,\omega)\le\mu(\omega),$$
hence the mixed index $\gamma(\sigma,\omega)$ can become crucial whenever $\gamma(\omega)<\mu(\omega)$ does hold true.
\end{remark}

The next result is analogous to Proposition \ref{gammamixedvsomegam} and showing that $\mu(\omega)$ is the upper value for our considerations.

\begin{proposition}\label{gammamixedvsomegamfct}
Let $\omega$ be a weight function (in the sense of Section \ref{section22}). Then we get
$$\sup\{r\in\RR_{>0}: \eqref{gamma1omegadescendent}\;\;\text{is satisfied}\}=\mu(\omega).$$
\end{proposition}

\demo{Proof}
If $\mu(\omega)=0$, then we have obviously equality. So let now $\mu(\omega)>0$.

As commented in $(iii)$ in Remark \ref{usefulremarks}, if \eqref{gamma1omegadescendent} with $r>0$ holds true, then $\omega^r$ has \hyperlink{omnq}{$(\omega_{\text{nq}})$} and consequently $\sup\{r\in\RR_{>0}: \eqref{gamma1omegadescendent}\;\;\text{is satisfied}\}\le\mu(\omega)$.\vspace{6pt}

Conversely, for each $0<r<\mu(\omega)$, we have that there exist some $t_0\ge 1$, $\varepsilon>0$ and $C\ge 1$ such that $\omega(t)\le C t^{1/r-\varepsilon}$ for all $t\ge t_0$. Thus for all $t\ge 0$ we get:
$$\int_1^{\infty}\frac{\omega(tu)}{u^{1+1/r}}du=A+\int_{t_0}^{\infty}\frac{\omega(tu)}{u^{1+1/r}}du\le A+Ct^{1/r-\varepsilon}\int_{t_0}^{\infty}\frac{1}{u^{1+\varepsilon}}du\le A_1+Bt^{1/r},$$
and so \eqref{gamma1omegadescendent} holds true with $\sigma(t):=t^{1/r}$.
\qed\enddemo

As it happens for weight sequences, the choice $\sigma(t):=t^{1/r}$ is not optimal. For each $0<r<\mu(\omega)$ we can define the weight $\kappa_{\omega^r}$, using here the notation
\begin{equation}\label{gammamixedvsomegamfctequ}
\kappa_{\omega}(t):=\int_1^{\infty}\frac{\omega(ty)}{y^2}dy=t\int_t^{\infty}\frac{\omega(y)}{y^2}dy.
\end{equation}
A weight $\omega$ satisfies \hyperlink{omsnq}{$(\omega_{\text{snq}})$} (or \hyperlink{gammarfct}{$(\omega_{\gamma_1})$}) if and only if $\kappa_{\omega}\hyperlink{sim}{\sim}\omega$. $\kappa_{\omega}$ is called a {\itshape heir} in \cite[Defintion 3.14]{whitneyextensionmixedweightfunction}. Thus \hyperlink{gammarfctmix}{$(\kappa^{1/r}_{\omega^r},\omega)_{\gamma_r}$} and more precisely $\kappa^{1/r}_{\omega^r}$ is optimal (i.e. minimal), up to a constant, among the weights satisfying requirement \eqref{gamma1omegadescendent}.


Note that $\kappa^{1/r}_{\omega^r}$ is a weight function (in the sense of Section \ref{section22}) and moreover \hyperlink{om1}{$(\omega_1)$}, \hyperlink{om3}{$(\omega_3)$}, \hyperlink{om4}{$(\omega_4)$} can be transferred from $\omega^r$ (so from $\omega$) to $\kappa_{\omega^r}$ and so to $\kappa^{1/r}_{\omega^r}$ (see Section \ref{section22}). The proof of \cite[Theorem 2.11 $(i)\Rightarrow(ii)$]{firstindexpaper} shows that $\kappa_{\omega^r}$ is concave. $(ii)$ in Remark \ref{usefulremarks} yields $\kappa^{1/r}_{\omega^r}(t)\ge\omega(t)$ for all $t\ge 0$.\vspace{6pt}

We are closing this section by establishing now the connection between $\mu(\omega)$ and $\mu(W^l)$, with $W^l\in\Omega$ and $\Omega$ denoting the matrix associated with $\omega$.

\begin{lemma}\label{Qvsomega}
Let $\omega\in\hyperlink{omset0}{\mathcal{W}_0}$ (i.e. a normalized weight function with \hyperlink{om3}{$(\omega_3)$} and \hyperlink{om4}{$(\omega_4)$}), let $\Omega=\{W^l: l>0\}$ be the matrix associated with $\omega$. Then we obtain
\begin{equation}\label{Qvsomegaequ}
\forall\;l>0:\;\;\;\mu(\omega)=\mu(\omega_{W^l})=\mu(W^l).
\end{equation}
In particular, if $\omega\equiv\omega_N$ for some given sequence $N\in\hyperlink{LCset}{\mathcal{LC}}$, then $\mu(\omega_N)=\mu(N)$.
\end{lemma}

The last statement yields that, if even $n\in\hyperlink{LCset}{\mathcal{LC}}$ holds true, then $\mu(\omega_N)=\mu(N)=\mu(n)+1=\mu(\omega_n)+1$.

\demo{Proof}
First, given $r>0$, by the formula on p. 7 in \cite{sectorialextensions1} we know that the matrix associated with the weight $\omega^r$ coincides with the set $\{V^{l,r}:=(W^{l/r})^{1/r}: l>0\}$, i.e taking the $r$-th root of the sequences belonging to $\Omega$ and making a re-parametrization of the index (in terms of $r$).

By \cite[Cor. 5.8 $(1)$]{compositionpaper}, which is based on the characterization given in \cite[Lemma 4.1]{Komatsu73} (see also \cite[Cor. 4.8]{testfunctioncharacterization}), applied to $\omega^r$ we know that $\omega^r$ satisfies \hyperlink{omnq}{$(\omega_{\text{nq}})$} if and only if some/each $V^{l,r}$ does have \hyperlink{mnq}{$(\text{nq})$} and so by the definition of the growth indices given in \eqref{orderofquasi} and \eqref{orderofquasifct} we are done.

If $\omega\equiv\omega_N$, then this follows immediately by recalling \eqref{omegaMspower} and the definitions of the growth indices.
\qed\enddemo

\begin{remark}
In \cite[Def. 3.3, Thm. 3.4]{Sanzflatultraholomorphic}, by using the connection between $\mu(N)$ and the so-called exponent of convergence of a nondecreasing sequence of positive real numbers tending to infinity (see \cite[Prop. 2.13]{Sanzflatultraholomorphic} and the references therein), for any $N\in\hyperlink{LCset}{\mathcal{LC}}$ we obtain $\mu(\widehat{N})=\mu(N)+1$, $\widehat{N}:=(p!N_p)_{p\in\NN}$.
\end{remark}

\section{Optimal flat functions in the mixed setting}\label{optimalflat}

\subsection{Construction of outer functions}\label{constructionouterfunction}

The aim of this paragraph is to obtain holomorphic functions in the right half-plane of $\CC$ whose growth is accurately controlled by two given weight functions $\omega$ and $\sigma$. Since in the forthcoming sections we want to treat the weight sequence and weight function case simultaneously we will transfer the general proofs from \cite[Section 6]{sectorialextensions} to a mixed setting.

First we translate \hyperlink{gammarfctmix}{$(\sigma,\omega)_{\gamma_1}$} into a property for $\sigma^{\iota}$ and $\omega^{\iota}$ (recall $\omega^{\iota}(t)=\omega(1/t)$). The next result should be compared with \cite[Lemma 6.1]{sectorialextensions} and \cite[Lemma 2.1.1]{Thilliezdivision}.

\begin{lemma}\label{functionh3tau}
Let $\omega$ and $\sigma$ be weight functions. Then, one has \hyperlink{gammarfctmix}{$(\sigma,\omega)_{\gamma_1}$} if, and only if,
\begin{equation*}
\exists\;C\ge 1\;\forall\;t>0:\;\;\;\int_0^1 -\omega^{\iota}(tu)du\ge-C(\sigma^{\iota}(t)+1).
\end{equation*}
\end{lemma}


The next statement recalls \cite[Lemma 6.2]{sectorialextensions}, see also \cite[Lemma 2.1.2]{Thilliezdivision}.

\begin{lemma}\label{Lemma212tau}
Let $\omega$ be a weight function. Then $\omega$ satisfies \hyperlink{omnq}{$(\omega_{\on{nq}})$} if and only if
\begin{equation}\label{Lemma212tauequ}
\int_{-\infty}^{+\infty}\frac{-\omega^{\iota}(|t|)}{1+t^2}dt>-\infty.
\end{equation}
In particular \eqref{Lemma212tauequ} holds true for each $\omega$ satisfying \hyperlink{gammarfctmix}{$(\sigma,\omega)_{\gamma_1}$} (and $\sigma$ some other possibly quasianalytic weight).
\end{lemma}


In the next step we are generalizing \cite[Lemma 6.3]{sectorialextensions} to a mixed setting, the idea of the construction in the proof is coming from \cite[Lemma 2.1.3]{Thilliezdivision}.

\begin{lemma}\label{Lemma213tau}
Let $\omega$ and $\sigma$ be two weight functions satisfying \hyperlink{gammarfctmix}{$(\sigma,\omega)_{\gamma_1}$}. Then for all $a>0$ there exists a function $F_a$ which is holomorphic on the right half-plane $H_1:=\{w\in\CC: \Re(w)>0\}\subseteq\CC$ and constants $A,B\ge 1$ (large) depending only on the weights $\omega$ and $\sigma$ such that
\begin{equation}\label{Lemma213equtau}
\forall\;w\in H_1:\;\;\;B^{-a}\exp(-2aB\sigma^{\iota}(B^{-1}\Re(w)))\le|F_a(w)|\le \exp\left(-\frac{a}{2}\omega^{\iota}(A|w|)\right).
\end{equation}
\end{lemma}

\demo{Proof}
We are following the lines of the proof of \cite[Lemma 6.3]{sectorialextensions}, see also \cite[Lemma 2.1.3]{Thilliezdivision} for the single weight sequence case. Since \hyperlink{gammarfctmix}{$(\sigma,\omega)_{\gamma_1}$} is valid, the weight $\omega$ has to satisfy \hyperlink{omnq}{$(\omega_{\on{nq}})$}. Hence for $w\in H_1$ we can put
$$F_a(w):=\exp\left(\frac{1}{\pi}\int_{-\infty}^{+\infty} \frac{-a\omega^{\iota}(|t|)}{1+t^2}\frac{itw-1}{it-w}dt\right);$$
Lemma \ref{Lemma212tau} implies immediately that $F_a$ is a holomorphic function in $H_1$. Since $F_a(w)=(F_1(w))^a$, we need only consider in the proof $a=1$ and put for simplicity $F:=F_1$.

For $w\in H_1$ write $w=u+iv$, hence $u>0$.
We have
$$\log(|F(w)|)=\frac{1}{\pi}\int_{\RR}-\omega^{\iota}(|t|) \frac{u}{(t-v)^2+u^2}dt=-\frac{1}{\pi}f\ast g_u(v),
$$
where $f(t):=\omega^{\iota}(|t|)$, $g_u(t):=u/(t^2+u^2)$. $f$ and $g_u$ are symmetrically nonincreasing functions, hence the convolution too, and as argued in \cite[Lemma 6.3]{sectorialextensions} this means that for $w\mapsto\log(|F(w)|)$ the minimum is attained on the positive real axis and we have for all $w\in H_1$:
$$\log(|F(w)|)\ge\log(|F(u)|)=\log(|F(\Re(w))|), \hspace{10pt}\log(|F(u)|)=\frac{1}{\pi}\int_{\RR}-\omega^{\iota}(|t|) \frac{u}{t^2+u^2}dt=-\frac{1}{\pi}f\ast g_u(0).$$
For the left-hand side in \eqref{Lemma213equtau} consider $K>0$ (small) and get
$$\pi\log(|F(u)|)=\int_{\RR}-\omega^{\iota}(|t|) \frac{u}{t^2+u^2}dt
=\int_{\{t:|t|\ge Ku\}}-\omega^{\iota}(|t|)\frac{u}{t^2+u^2}dt+\int_{\{t:|t|\le Ku\}}-\omega^{\iota}(|t|)\frac{u}{t^2+u^2}dt.$$
The first integral is estimated by
\begin{align*}
\int_{\{|t|\ge Ku\}}-\omega^{\iota}(|t|)\frac{u}{t^2+u^2}dt
&\ge-\omega^{\iota}(Ku)\int_{\{|t|\ge Ku\}}\frac{u}{t^2+u^2}dt
=-\omega^{\iota}(Ku)(\pi-2\arctan(K))
\\&
\ge-\sigma^{\iota}(Ku)D(\pi-2\arctan(K))-D(\pi-2\arctan(K)),
\end{align*}
since $t\mapsto -\omega^{\iota}(t)$ is nondecreasing and since by \hyperlink{gammarfctmix}{$(\sigma,\omega)_{\gamma_1}$} we have $\sigma\hyperlink{ompreceq}{\preceq}\omega$. Thus $-\omega^{\iota}(t)\ge-D(\sigma^{\iota}(t)+1)$ for some $D\ge 1$ and all $t>0$ follows.\vspace{6pt}

For the second integral we get
$$\int_{\{t: |t|\le Ku\}}-\omega^{\iota}(|t|)\frac{u}{t^2+u^2}dt
=\int_{\{s: |s|\le1\}}-\omega^{\iota}(Ku|s|)\frac{K}{K^2s^2+1}ds
\ge K\int_{\{s: |s|\le 1\}}-\omega^{\iota}(Ku|s|)ds,$$
since $-\omega^{\iota}(Ku|s|)\le 0$ holds for any $K,u>0$ and $|s|\le 1$. Let $C\ge 1$ be the constant appearing in Lemma \ref{functionh3tau}, then
$$K\int_{\{s:|s|\le 1\}}-\omega^{\iota}(Ku|s|)ds\ge2KC(-\sigma^{\iota}(Ku)-1).$$
Now follow the lines of the proof in \cite[Lemma 6.3]{sectorialextensions} with $-\tau^{\iota}$ replaced by $-\sigma^{\iota}$ which is also nondecreasing. Thus the first estimate in \eqref{Lemma213equtau} is shown.\vspace{6pt}

The second estimate follows without any change like in \cite[Lemma 6.3]{sectorialextensions}.

\qed\enddemo

\subsection{Construction of optimal sectorially flat functions in the mixed setting}\label{sectoriallyflatfunction}
Using the results from Section \ref{constructionouterfunction} the aim is now to transfer \cite[Theorem 2.3.1]{Thilliezdivision} and \cite[Theorem 6.7]{sectorialextensions} (see also \cite[Theorem 5.6]{sectorialextensions1}) simultaneously to a mixed setting.

\begin{theorem}\label{Theorem231mixed}
\begin{itemize}
\item[$(I)$] Let $M,N\in\hyperlink{LCset}{\mathcal{LC}}$ be given such that $\mu_p\le\nu_p$ and $\gamma(M,N)>0$ holds true.

    Then for any $0<\gamma<\gamma(M,N)$ there exist constants $K_1,K_2,K_3,K_4>0$ depending only on $M$, $N$ and $\gamma$ such that for all $a>0$ there exists a function $G_a$ holomorphic on $S_{\gamma}$ and satisfying
\begin{equation}\label{Theorem231sequequ}
\forall\;\xi\in S_{\gamma}:\;\;K_1^{-a}h_M(K_1|\xi|)^{2aK_2}\le|G_a(\xi)|\le h_N(K_3|\xi|)^{\frac{aK_4}{2}}.
\end{equation}
Moreover, if $N$ has in addition \hyperlink{mg}{$(\on{mg})$}, then $G_a\in\mathcal{A}_{\{\widehat{N}\}}(S_{\gamma})$ with $\widehat{N}:=(p!N_p)_{p\in\NN}$ (and $G_a$ is flat at $0$).

If $M$ has in addition \hyperlink{mg}{$(\on{mg})$}, then there exists $K_5>0$ depending also on given $a>0$ such that
\begin{equation}\label{Theorem231sequequ1}
\forall\;\xi\in S_{\gamma}:\;\;\;|G_a(\xi)|\ge K_1^{-a}h_M(K_5|\xi|).
\end{equation}

\item[$(II)$] Let $\omega$ and $\sigma$ be weight functions such that $\gamma(\sigma,\omega)>0$ holds true. Then for any $0<\gamma<\gamma(\sigma,\omega)$ there exist constants $K_1,K_2,K_3>0$ depending only on $\sigma$, $\omega$ and $\gamma$ such that for all $a>0$ there exists a function $G_a$ holomorphic on $S_{\gamma}$ and satisfying
\begin{equation}\label{Theorem231fctequ}
\forall\;\xi\in S_{\gamma}:\;\;\;K_1^{-a}\exp(-2a\sigma^{\iota}(K_2|\xi|))\le|G_a(\xi)|\le\exp\left(-\frac{a}{2}\omega^{\iota}(K_3|\xi|)\right).
\end{equation}
Moreover, if $\omega$ is normalized and satisfies \hyperlink{om3}{$(\omega_3)$}, then $G_a\in\mathcal{A}_{\{\widehat{\Omega}\}}(S_{\gamma})$ (and $G_a$ is flat at $0$), more precisely
\begin{equation}\label{Theorem231fctequ0}
\forall\;\xi\in S_{\gamma}:\;\;\;|G_a(\xi)|\le h_{W^{2/a}}(K'_2|\xi|),
\end{equation}
where $\Omega=\{W^x: x>0\}$ shall denote the matrix associated with $\omega$ and $\widehat{\Omega}:=\{\widehat{W}^x=(p!W^x_p)_{p\in\NN}: x>0\}$.

If $\sigma\in\hyperlink{omset0}{\mathcal{W}_0}$, then there exist an index $x>0$ and a constant $K_4>0$ depending also on $a$ such that
\begin{equation}\label{Theorem231fctequ1}
\forall\;\xi\in S_{\gamma}:\;\;\;|G_a(\xi)|\ge K_4 h_{S^x}(K_2|\xi|),
\end{equation}
where $S^x\in\Sigma$, $\Sigma$ the matrix associated with $\sigma$.
\end{itemize}
\end{theorem}

\demo{Proof}
We will give some more details for the proof of $(I)$ (following the lines of \cite[Theorem 2.3.1]{Thilliezdivision}).

$(I)$ Let $a>0$ be arbitrary. Take $s,\delta>0$ such that $\gamma<\delta<\gamma(M,N)$ and $s\delta<1<s\gamma(M,N)$.

By \eqref{growthindexpower} we get $s\gamma(M,N)>1\Leftrightarrow\gamma(M^s,N^s)>1$, hence, by Lemma \ref{mglemma1} applied to $M^s$ and $N^s$, we get $\gamma(\omega_{M^s},\omega_{N^s})>1$, too.

So we can use Lemma \ref{Lemma213tau} for $\omega_{M^s}\equiv\sigma$ and $\omega_{N^s}\equiv\omega$ and obtain a function $F_a$ holomorphic on the right half-plane and satisfying
\begin{equation}\label{Theorem231equ1}
\forall\;w\in H_1:\;\;\;B^{-a}\exp(-2aB\omega_{M^s}^{\iota}(B^{-1}\Re(w)))\le|F_a(w)|\le \exp\left(-\frac{a}{2}\omega_{N^s}^{\iota}(A|w|)\right).
\end{equation}
Then put
$$G_a(\xi)=F_a(\xi^s),\;\;\;\;\xi\in S_{\delta}.$$
Note that, as $s\delta<1$, the ramification $\xi\mapsto\xi^s$ maps holomorphically $S_{\delta}$ into $S_{\delta s}\subseteq S_1=H_1$, and so $G_a$ is well-defined. We show that the restriction of $G_a$ to $S_{\gamma}\subseteq S_{\delta}$ satisfies the desired properties by proving that \eqref{Theorem231sequequ} holds indeed on the whole $S_{\delta}$.

First we consider the lower estimate.

Let $\xi\in S_{\delta}$ be given, then $\Re(\xi^s)\ge\cos(s\delta\pi/2)|\xi|^s$ (since $s\delta\pi/2<\pi/2$). If $B\ge 1$ denotes the constant coming from the left-hand side in \eqref{Lemma213equtau} applied to the weight $\omega_{M^s}$ (or see \eqref{Theorem231equ1}), then
\begin{align*}
|G_a(\xi)|&=|F_a(\xi^s)|\ge B^{-a}\exp(-2aB\omega_{M^s}^{\iota}(B^{-1}(\Re(\xi^s))))\\
&\ge B^{-a}\exp(-2aB\omega_{M^s}^{\iota}(B^{-1}\cos(s\delta\pi/2)|\xi|^s))
\\&
=B^{-a}\exp(-2aB\omega_{M^s}^{\iota}((B_1|\xi|)^s))
=B^{-a}\exp(-2aBs\omega_M^{\iota}(B_1|\xi|))=B^{-a}h_M(B_1|\xi|)^{2aBs},
\end{align*}
where we have put $B_1:=(B^{-1}\cos(s\delta\pi/2))^{1/s}$. We have used $\omega_{M^s}^{\iota}(t^s)=\omega_{M^s}(1/t^s)=s\omega_M(1/t)=s\omega_M^{\iota}(t)$ for all $t,s>0$, see \eqref{omegaMspower}, and finally \eqref{functionhequ2}.\vspace{6pt}

Now we consider the right-hand side in \eqref{Lemma213equtau} respectively in \eqref{Theorem231equ1} and proceed as before. Let $A$ be the constant coming from the right-hand side of \eqref{Lemma213equtau} applied to $\omega_{N^s}$, so
\begin{align*}
|G_a(\xi)|&=|F_a(\xi^s)|\le \exp\left(-\frac{a}{2}\omega_{N^s}^{\iota}(A|\xi|^s)\right)
=\exp\left(-\frac{a}{2}\omega_{N^s}^{\iota}((A^{1/s}|\xi|)^s)\right)
\\&
=\exp\left(-\frac{a}{2}s\omega_{N}^{\iota}(A^{1/s}|\xi|)\right)=h_N(A^{1/s}|\xi|)^{\frac{sa}{2}},
\end{align*}
and \eqref{Theorem231sequequ} has been proved for every $\xi\in S_{\delta}$.

In order to show $G_a\in\mathcal{A}_{\{\widehat{N}\}}(S_{\gamma})$ we put in the estimate above $A_1:=A^{1/s}$ and see
$$\exp\left(-\frac{a}{2}s\omega_N^{\iota}(A_1|\xi|)\right)=\exp\left(-\frac{a}{2}s\omega_N(1/(A_1|\xi|))\right)=h_N(A_1|\xi|)^{sa/2}.$$
If we can show
\begin{equation}\label{equaGsubaflat}
\exists\;A_2\ge 1\;\forall \xi\in S_{\delta}:\ \ |G_a(\xi)|\le h_{N}(A_2|\xi|),
\end{equation}
then by applying \cite[Lemma 6.4 $(i.1)$]{sectorialextensions} we see that $G_a\in\mathcal{A}_{\{\widehat{N}\}}(S_{\gamma})$ (and it is a flat function at $0$). Since $h_N\le 1$, \eqref{equaGsubaflat} holds true whenever $\frac{sa}{2}\ge 1\Leftrightarrow sa\ge 2$. But in general we have to use \hyperlink{mg}{$(\text{mg})$} for $N$ and iterate \eqref{mgsquare} (applied for $N$) $l$-times, $l\in\NN$ chosen minimal to ensure $\frac{sa}{2}\ge\frac{1}{2^l}$.

The proof of \eqref{Theorem231sequequ1} follows analogously by iterating \hyperlink{mg}{$(\text{mg})$} for $M$ (if necessary) in order to get rid of the exponent $2aK_2$.\vspace{6pt}

$(II)$ Again let $a>0$ be arbitrary but fixed, take $s,\delta>0$ such that $\gamma<\delta<\gamma(\sigma,\omega)$, $s\delta<1<s\gamma(\sigma,\omega)=\gamma(\sigma^{1/s},\omega^{1/s})$ and put
$$G_a(\xi)=F_a(\xi^s),\;\;\;\;\xi\in S_{\delta},$$
where $F_a$ is the function constructed in Lemma \ref{Lemma213tau}. We apply Lemma \ref{Lemma213tau} to $\sigma^{1/s}\equiv\sigma$, $\omega^{1/s}\equiv\omega$. Hence \eqref{Theorem231fctequ} holds true by following the lines in the proof of \cite[Theorem 6.7]{sectorialextensions} and replacing $\tau$ by $\omega$ for the right-hand side respectively $\tau$ by $\sigma$ for the left-hand side.\vspace{6pt}

The remaining statements, in particular the estimate \eqref{Theorem231fctequ1}, follow analogously as in the proof of \cite[Theorem 6.7]{sectorialextensions}, replacing $\tau$ by $\omega$ or $\sigma$ in the arguments.

\qed\enddemo

\section{Right inverses for the asymptotic Borel map in ultraholomorphic classes in sectors}\label{sectRightInver1var}

The aim of this section is to obtain an extension result in the ultraholomorphic classes considered in a mixed setting for both the weight sequence and the weight function approach following the proofs and techniques in \cite[Section 7]{sectorialextensions}. The existence of the optimal flat functions $G_a$ obtained in Theorem~\ref{Theorem231mixed} will be the main ingredient in the proof which is inspired by the same technique as in previous works of A. Lastra, S. Malek and the second author~\cite{LastraMalekSanzContinuousRightLaplace, Sanzsummability} in the single weight sequence approach. Although for the general construction the weight functions $\sigma$ and $\omega$ need not be normalized, we are interested in working with the weight matrices associated with them, which will be standard log-convex if we ask for normalization and \hyperlink{om3}{$(\omega_3)$} to hold.

Note that any weight function may be substituted by a normalized equivalent one (e.g. see \cite[Remark 1.2 $(b)$]{BonetBraunMeiseTaylorWhitneyextension}) and equivalence preserves the property \hyperlink{om3}{$(\omega_3)$}, so it is no restriction to ask for normalization from the very beginning.

\begin{remark}\label{omega1loss}
An important difference to the complete approach in \cite{sectorialextensions} is, see also the comments given in the introduction in Section 7 there, that condition $\gamma(\omega)>0$ and which amounts to \hyperlink{om1}{$(\omega_1)$} as shown in \cite[Lemma 4.2]{sectorialextensions} will not be valid in general anymore in the mixed situation. In the following we are only requiring $\gamma(\sigma,\omega)>0$ and recall that $\gamma(\sigma,\omega)\ge\gamma(\omega)$ as shown above. An explicit example of this situation, having $\gamma(\sigma,\omega)>0$ (as large as desired) and $\gamma(\omega)=\gamma(\sigma)=0$ will be provided in the Appendix \ref{examplesfromlangenbruch} below. We are able to treat this situation by recognizing that in \cite{sectorialextensions} we have worked in a very general framework for weight functions and the assumption $\gamma(\omega)>0$ can be replaced by $\gamma(\sigma,\omega)>0$ without causing problems.\vspace{6pt}

Recall that \hyperlink{om1}{$(\omega_1)$} is standard in the ultradifferentiable setting and thus our techniques make it possible to treat ''exotic'' weight function situations as well. Moreover \hyperlink{om1}{$(\omega_1)$} has also been used to have that the class defined by $\omega$ admits a representation by using the associated weight matrix $\Omega$, see Section \ref{weightmatrixfromfunction}. Thus the warranty that the ultraholomorphic (and also the ultradifferentiable) spaces associated with $\omega$ and its corresponding weight matrix $\Omega$ coincide is not clear anymore, see the comments preceding~\eqref{equaEqualitySpacesWeightFunctionMatrix}.

Therefore the main and most general ultraholomorphic extension result Theorem \ref{theoExtensionOperatorsMatrixmixed} deals with a mixed situation between classes defined by (associated) weight matrices. If one imposes \hyperlink{om1}{$(\omega_1)$} on the weights one is able to prove a mixed version of classes defined by weight functions, see Corollary \ref{coroExtensionOperatorsWeights}. Finally, in Theorem \ref{theoExtensionOperatorssequencemixed} we will treat the mixed weight sequence case as well.
\end{remark}

\subsection{Preliminaries}
We start with recalling \cite[Lemma 7.1]{sectorialextensions}.

\begin{lemma}\label{lemmaKernelsmixed}
Let $\sigma$ and $\omega$ be normalized weight functions with $\gamma(\sigma,\omega)>0$ and such that $\omega$ satisfies
\hyperlink{om3}{$(\omega_3)$}. Let $\Omega=\{W^x: x>0\}$ be the weight matrix associated with $\omega$, $0<\gamma<\gamma(\sigma,\omega)$, and for $a>0$ let $G_a$ be the function constructed in Theorem~\ref{Theorem231mixed}.
Let us define the function $e_{a}:S_{\gamma}\to\CC$  by
$$
e_a(z):=z\,G_a(1/z),\quad z\in S_{\gamma}.
$$
The function $e_{a}$ enjoys the following properties:
\begin{itemize}
\item[(i)] $z^{-1}e_{a}(z)$ is uniformly integrable at the origin, it is to say, for any $t_0>0$ we have

    $$\int_{0}^{t_0}t^{-1}\sup_{|\tau|<\gamma\pi/2}|e_{a}(te^{i\tau})|dt<\infty.$$
\item[(ii)] There exist constants $K>0$, independent from $a$, and $C>0$, depending on $a$, such that
\begin{equation}\label{equaBounds.for.kernel}
|e_{a}(z)|\le Ch_{W^{4/a}}\left(\frac{K}{|z|}\right),\qquad z\in S_{\gamma}.
\end{equation}
\item[(iii)] For $\xi\in\RR$, $\xi>0$, the values of $e_{a}(\xi)$ are positive real.
\end{itemize}
\end{lemma}

\demo{Proof}
The proof is completely the same as for \cite[Lemma 7.1]{sectorialextensions}, for $(i)$ we apply the right-hand side in \eqref{Theorem231fctequ}, for $(ii)$ we use \eqref{Theorem231fctequ0} and \eqref{functionh2equ1} together with the definition given in \eqref{functionhequ1}.




\qed\enddemo

Analogously as in \cite[Definition 7.2]{sectorialextensions} we introduce now the moment function associated with $e_a$.

\begin{definition}
We define the \textit{moment function} associated with the function $e_a$ (introduced in the previous Lemma) as
$$m_a(\lambda):=\int_{0}^{\infty}t^{\lambda-1}e_{a}(t)dt=
\int_{0}^{\infty}t^{\lambda}G_{a}(1/t)dt.$$
\end{definition}

From Lemma~\ref{lemmaKernelsmixed} and the definition of $h_{W^x}$ in \eqref{functionhequ1} we see that for every $p\in\NN$,
$$
|e_{a}(z)|\le C\frac{K^pW^{4/a}_p}{|z|^p},\qquad z\in S_{\gamma}.
$$
So, we easily deduce that the function $m_a$ is well defined and continuous in $\{\lambda: \Re(\lambda)\ge0\}$, and holomorphic in $\{\lambda: \Re(\lambda)>0\}$. Moreover, $m_a(\xi)$ is positive for every $\xi\ge0$, and the sequence $(m_a(p))_{p\in\NN}$ is called the \textit{sequence of moments} of $e_{a}$.

The next result is generalizing \cite[Proposition 7.3]{sectorialextensions}, which is similar to Proposition 3.6 in~\cite{LastraMalekSanzContinuousRightLaplace}, to a mixed setting.

\begin{proposition}\label{propmequivmmixed}
Let $\sigma$ and $\omega$ be normalized weight functions with $\gamma(\sigma,\omega)>0$ and such that both weights satisfy
\hyperlink{om3}{$(\omega_3)$}. Let $\Sigma=\{S^x: x>0\}$ and $\Omega=\{W^x: x>0\}$ be the weight matrices associated with $\sigma$ and $\omega$ respectively, and for $0<\gamma<\gamma(\sigma,\omega)$ and $a>0$ let $G_a,e_a,m_a$ be the functions previously constructed. Then, there exist constants $C_1,C_2>0$, both depending on $a$, such that for every $p\in\NN$ one has
\begin{equation}
  \label{equaIneqMomentsMatrixmixed}
  C_1\left(\frac{K_2}{2}\right)^pS_p^{1/(2a)}\le m_a(p)\le C_2K_3^pW_{p}^{4/a},
\end{equation}
where $K_2$ and $K_3$ are the constants, not depending on $a$, appearing in \eqref{Theorem231fctequ}.
\end{proposition}

\demo{Proof}
The proof follows the lines as in \cite[Proposition 7.3]{sectorialextensions} (based on the arguments by O. Blasco in~\cite{Blasco}). For the second estimate in \eqref{equaIneqMomentsMatrixmixed} we use the second inequality in \eqref{Theorem231fctequ} (and here also \eqref{newmoderategrowth} is used); the first estimate in \eqref{equaIneqMomentsMatrixmixed} follows by applying the first inequality in \eqref{Theorem231fctequ}.

\qed\enddemo


\subsection{Main extension results}
Now we are able to formulate and proof the generalization of the main extension result \cite[Theorem 7.4]{sectorialextensions}.

\begin{theorem}\label{theoExtensionOperatorsMatrixmixed}
Let $\sigma$ and $\omega$ be normalized weight functions with $\gamma(\sigma,\omega)>0$ and such that both weights satisfy
\hyperlink{om3}{$(\omega_3)$} and $0<\gamma<\gamma(\sigma,\omega)$. Moreover we denote by $\Sigma=\{S^x: x>0\}$ and $\Omega=\{W^x: x>0\}$ the weight matrices associated with $\sigma$ and $\omega$ respectively and consider the matrices $\widehat{\Sigma}=\{\widehat{S}^x:x>0\}$ and $\widehat{\Omega}=\{\widehat{W}^x:x>0\}$ where $\widehat{S}^x:=(p!S^x_p)_{p\in\NN}$ and $\widehat{W}^x:=(p!W^x_p)_{p\in\NN}$.

Then, there exists a constant $k_0>0$ such that for every $x>0$ and every $h>0$, one can construct a linear and continuous map
$$
E^{\sigma,\omega}_h:\;\;\lambda\in\Lambda_{\widehat{S}^x,h}\mapsto f_{\lambda}\in \mathcal{A}_{\widehat{W}^{8x},k_0h}(S_{\gamma})
$$
such that for every $\lambda\in\Lambda_{\widehat{S}^x,h}$ one has $\mathcal{B}\circ E^{\sigma,\omega}_h(\lambda)=\mathcal{B}(f_{\lambda})=\lambda$. Consequently we obtain the inclusion $\mathcal{B}(\mathcal{A}_{\{\widehat{\Omega}\}}(S_{\gamma}))\supseteq\Lambda_{\{\widehat{\Sigma}\}}$.
\end{theorem}

\demo{Proof}
Fix $\delta>0$ such that $\gamma<\delta<\gamma(\sigma,\omega)$. Given $\lambda=(\lambda_p)_{p\in\NN}\in\Lambda_{\widehat{S}^x,h}$, we have
\begin{equation}\label{equaNormlambdamixed}
|\lambda_p|\le |\lambda|_{\widehat{S}^x,h} h^{p}p!S_p^x,\quad p\in\NN_0.
\end{equation}
We choose $a=1/(2x)$, and consider the function $G_a$, defined in $S_{\delta}$, obtained in Theorem~\ref{Theorem231mixed} for such value of $a$, and the corresponding functions $e_a$ and $m_a$ previously defined.
Next, we consider the formal power series
$$
\widehat{f}_{\lambda}:=\sum_{p=0}^\infty\frac{\lambda_p}{p!}z^p
$$
and its formal (Borel-like) transform
$$
\widehat{\mathcal{B}}_a\widehat{f}_{\lambda}:= \sum_{p=0}^\infty\frac{\lambda_p}{p!m_a(p)}z^p.
$$
By the choice of $a$, \eqref{equaNormlambdamixed} and the first part of the inequalities in~\eqref{equaIneqMomentsMatrixmixed}, we deduce that
\begin{equation}
  \label{equaSizeCoefficientsBorelTransformmixed}
  \left|\frac{\lambda_p}{p!m_a(p)}\right|\le
  \frac{|\lambda|_{\widehat{S}^x,h} h^{p}p!S_p^x}
  {C_1\left(K_2/2\right)^p p!S_p^x}=
  \frac{|\lambda|_{\widehat{S}^x,h}}{C_1}
  \left(\frac{2h}{K_2}\right)^p,
\end{equation}
and so the series $\widehat{\mathcal{B}}_a\widehat{f}_{\lambda}$ converges in the disc of center 0 and radius $K_2/(2h)$ (not depending on $\lambda$), where it defines a holomorphic function $g_\lambda$. We set $R_0:=K_2/(4h)$, and define
\begin{equation*}
f_{\lambda}(z):=\int_{0}^{R_0}e_{a}\left(\frac{u}{z}\right)g_{\lambda}(u)\frac{du}{u},\qquad z\in S_{\delta}.
\end{equation*}
By virtue of Leibniz's theorem on analyticity of parametric integrals, $f_\lambda$ is holomorphic in $S_{\delta}$.

Our next aim is to obtain suitable estimates for the difference between $f$ and the partial sums of the series $\widehat{f}_{\lambda}$. As in the non-mixed setting, for given  $N\in\NN_0$ and $z\in S_{\delta}$ we have
\begin{align*}
f_{\lambda}(z)-\sum_{p=0}^{N-1}\lambda_p\frac{z^p}{p!}= \int_{0}^{R_0}e_{a}\left(\frac{u}{z}\right) \sum_{p=0}^{\infty}\frac{\lambda_{p}}{m_a(p)}\frac{u^p}{p!} \frac{du}{u} -\sum_{p=0}^{N-1}\frac{\lambda_p}{m_a(p)} \int_{0}^{\infty}u^{p-1}e_{a}(u)du\frac{z^p}{p!}.
\end{align*}
In the second integral we make the change of variable $v=zu$, what results in a rotation of the line of integration. By the estimate~(\ref{equaBounds.for.kernel}), one may use Cauchy's residue theorem in order to obtain that
$$
z^p\int_{0}^{\infty}u^{p-1}e_{a}(u)du= \int_{0}^{\infty}v^{p}e_{a}\left(\frac{v}{z}\right)\frac{dv}{v},
$$
which allows us to write the preceding difference as
\begin{multline*}
\int_{0}^{R_0}e_{a}\left(\frac{u}{z}\right) \sum_{p=0}^{\infty}\frac{\lambda_{p}}{m_a(p)} \frac{u^p}{p!}\frac{du}{u} -\sum_{p=0}^{N-1}\frac{\lambda_p}{m_a(p)} \int_{0}^{\infty}u^{p}e_{a}\left(\frac{u}{z}\right)\frac{du}{u} \frac{1}{p!}\\
=\int_{0}^{R_0}e_{a}\left(\frac{u}{z}\right) \sum_{p=N}^{\infty}\frac{\lambda_{p}}{m_a(p)}\frac{u^p}{p!} \frac{du}{u} -\int_{R_0}^{\infty}e_{a}\left(\frac{u}{z}\right) \sum_{p=0}^{N-1}\frac{\lambda_p}{m_a(p)} \frac{u^{p}}{p!}\frac{du}{u}.
\end{multline*}
Then, we have
\begin{equation}\label{equaBoundsAsympExpanflambda}
  \left|f_{\lambda}(z)-\sum_{p=0}^{N-1}\lambda_p\frac{z^p}{p!}\right|\le |f_{1}(z)|+|f_2(z)|,
\end{equation}
where
$$f_{1}(z)=\int_{0}^{R_0}e_{a}\left(\frac{u}{z}\right) \sum_{p=N}^{\infty}\frac{\lambda_{p}}{m_a(p)} \frac{u^p}{p!}\frac{du}{u},\quad
f_{2}(z)=\int_{R_0}^{\infty}e_{a}\left(\frac{u}{z}\right) \sum_{p=0}^{N-1}\frac{\lambda_p}{m_a(p)} \frac{u^{p}}{p!}\frac{du}{u}.$$

From~\eqref{equaSizeCoefficientsBorelTransformmixed} we deduce that
\begin{align}
|f_{1}(z)|&\le \frac{|\lambda|_{\widehat{S}^x,h}}{C_1}
  \int_{0}^{R_0}\left|e_{a}\left(\frac{u}{z}\right)\right| \sum_{p=N}^{\infty}\Big(\frac{2hu}{K_2}\Big)^p\frac{du}{u}
  =\frac{|\lambda|_{\widehat{S}^x,h}}{C_1} \Big(\frac{2h}{K_2}\Big)^N
  \int_{0}^{R_0}\left|e_{a}\left(\frac{u}{z}\right)\right| \frac{u^N}{1-\frac{2hu}{K_2}}\frac{du}{u}\nonumber\\
  &\le \frac{2|\lambda|_{\widehat{S}^x,h}}{C_1} \Big(\frac{2h}{K_2}\Big)^N
  \int_{0}^{R_0}\left|e_{a}\left(\frac{u}{z}\right)\right| u^{N-1}\,du,\label{equaBoundsf-1}
\end{align}
where in the last step we have used that $0<u<R_0=K_2/(4h)$ we have $1-2hu/K_2>1/2$.
In order to estimate $f_{2}(z)$, observe that for $u\ge R_0$ and $0\le p\le N-1$ we always have $u^p\le R_0^pu^N/R_0^N$, and so, using again~\eqref{equaSizeCoefficientsBorelTransformmixed} and the value of $R_0$, we may write
$$\left|\sum_{p=0}^{N-1}\frac{\lambda_pu^p}{p!m_a(p)}\right|\le \frac{|\lambda|_{\widehat{S}^x,h}}{C_1}\frac{u^N}{R_0^N} \sum_{p=0}^{N-1}\Big(\frac{2h}{K_2}\Big)^p R_0^p
\le \frac{2|\lambda|_{\widehat{S}^x,h}}{C_1} \Big(\frac{4h}{K_2}\Big)^N u^N.
$$
Then, we deduce that
\begin{equation}\label{equaBoundsf-2}
  |f_2(z)|\le \frac{2|\lambda|_{\widehat{S}^x,h}}{C_1} \Big(\frac{4h}{K_2}\Big)^N \int_{R_0}^{\infty}\left|e_{a}\left(\frac{u}{z}\right)\right| u^{N-1}du.
\end{equation}
In order to conclude, it suffices then to obtain estimates for $\int_{0}^{\infty}|e_{a}(u/z)|u^{N-1}du$. For this, note first that, by the estimates in~\eqref{Theorem231fctequ},
\begin{align*}
  \int_{0}^{\infty}\left|e_{a}\left(\frac{u}{z}\right)\right| u^{N-1}du&=
  \int_{0}^{\infty}\frac{u}{|z|} \left|G_{a}\left(\frac{z}{u}\right)\right| u^{N-1}du\\
  &\le \int_0^\infty \frac{u^N}{|z|}\exp\left(-\frac{a}{2}\omega\left(\frac{u}{K_3|z|}\right)\right)\,du=
  |z|^N\int_0^\infty t^N\exp\left(-\frac{a}{2}\omega\left(\frac{t}{K_3}\right)\right)\,dt.
\end{align*}
Now, we can follow the first part of the proof of \cite[Proposition 7.3]{sectorialextensions} to obtain that
\begin{equation}
  \label{equaBoundsIntegralsKerneltimespower}
  \int_{0}^{\infty}\left|e_{a}\left(\frac{u}{z}\right)\right| u^{N-1}du\le C_2K_3^NW_N^{4/a}|z|^N=C_2K_3^NW_N^{8x}|z|^N.
\end{equation}
Gathering~\eqref{equaBoundsAsympExpanflambda}, \eqref{equaBoundsf-1}, \eqref{equaBoundsf-2} and \eqref{equaBoundsIntegralsKerneltimespower}, we get
\begin{equation}
  \label{equaAsExflambdafinal}
  \left|f_{\lambda}(z)-\sum_{p=0}^{N-1}\lambda_p\frac{z^p}{p!}\right|\le \frac{2C_2|\lambda|_{\widehat{S}^x,h}}{C_1} \left(\frac{4hK_3}{K_2}\right)^NW_N^{8x}|z|^N.
\end{equation}
A straightforward application of Cauchy's integral formula yields that there exists a constant $r$, depending only on $\gamma$ and $\delta$, such that whenever $z$ is restricted to belong to $S_{\gamma}$, one has that
for every $p\in\NN$,
\begin{equation*}
|f^{(p)}(z)|\le \frac{2C_2|\lambda|_{\widehat{S}^x,h}}{C_1} \left(\frac{4hK_3r}{K_2}\right)^p p!W_p^{8x}.
\end{equation*}
So, putting $k_0:=\frac{4K_3r}{K_2}$ (independent from $x$ and $h$), we see that
$f_{\lambda}\in \mathcal{A}_{\widehat{W}^{8x},k_0h}(S_{\gamma})$ and $\|f_{\lambda}\|_{\widehat{W}^{8x},k_0h}\le
\frac{2C_2}{C_1}|\lambda|_{\widehat{S}^x,h}$. Since the map sending $\lambda$ to $f_{\lambda}$ is clearly linear, this last inequality implies that the map is also continuous from $\Lambda_{\widehat{S}^x,h}$ into
$\mathcal{A}_{\widehat{W}^{8x},k_0h}(S_{\gamma})$. Finally, from \eqref{equaAsExflambdafinal} one may easily deduce that
$\mathcal{B}(f_{\lambda})=\lambda$, and we conclude.
\qed\enddemo


\begin{remark}\label{notpossible}
We point out that by checking the proofs of the main extension results \cite[Theorem 3.2.1]{Thilliezdivision} and \cite[Theorem 6.12]{sectorialextensions1} it seems to be not possible to transfer the techniques there to a mixed setting. In these proofs the existence of a continuous linear extension operator coming from the ultradifferentiable extension theorems is used and it is sending the sequence $\lambda$ to the function $g_{\lambda}$ (e.g. see \cite[(6.12)]{sectorialextensions1}) which belongs in a mixed situation to a space defined by a different weight sequence or weight function.

Consequently, in the estimates of part $(ii)$ in the proofs of \cite[Theorem 3.2.1]{Thilliezdivision} and \cite[Theorem 6.12]{sectorialextensions1}, one is not able to get rid of the quotients of associated functions $h_M(\cdot)$. This should be compared with the estimate \eqref{equaSizeCoefficientsBorelTransformmixed} where we do not have to change the weight structure and the proofs of the extension theorems in \cite{ChaumatChollet94} and \cite{whitneyextensionmixedweightfunction}.
\end{remark}

In order to prove the generalization of \cite[Corollary 7.6]{sectorialextensions} first we have to recall \cite[Theorem 5.3]{sectorialextensions}.

\begin{theorem}\label{theoEqualSpacesMatrixwithfactorialWeight}
Let $\omega\in\hyperlink{omset1}{\mathcal{W}}$, i.e. a normalized weight satisfying \hyperlink{om3}{$(\omega_3)$}, \hyperlink{om4}{$(\omega_4)$} and \hyperlink{om1}{$(\omega_1)$}, $\Omega=\{W^x: x>0\}$ be the weight matrix associated with $\omega$ and consider $\widehat{\Omega}=\{\widehat{W}^x:x>0\}$ where $\widehat{W}^x=(p!W^x_p)_{p\in\NN}$. Then the following identities hold as locally convex vector spaces for all sector $S$ and for all $x>0$:
$$\mathcal{A}_{\{\widehat{\Omega}\}}(S)= \mathcal{A}_{\{\omega_{\widehat{W}^x}\}}(S),$$
and the same equalities are valid for the corresponding sequence classes $\Lambda$.
So, $\mathcal{A}_{\{\widehat{\Omega}\}}(S)$ coincides with the space $\mathcal{A}_{\{\tau\}}(S)$ associated with $\tau=\omega_{\widehat{W}^x}\in\hyperlink{omset}{\mathcal{W}}$.

Finally one gets $\gamma(\tau)=\gamma(\omega)+1>1$.
\end{theorem}

Now we concentrate on the generalization of \cite[Corollary 7.6]{sectorialextensions}.

\begin{corollary}\label{coroExtensionOperatorsWeights}
Let $\sigma, \omega\in\hyperlink{omset1}{\mathcal{W}}$ be given, so that $\gamma(\sigma,\omega)\ge\gamma(\omega)>0$, and let $0<\gamma<\gamma(\sigma,\omega)$ and $\Sigma=\{S^x: x>0\}$ and $\Omega=\{W^x: x>0\}$ be the weight matrices associated with $\sigma$ and $\omega$ respectively and consider the matrices $\widehat{\Sigma}=\{\widehat{S}^x:x>0\}$ and $\widehat{\Omega}=\{\widehat{W}^x:x>0\}$ where $\widehat{S}^x=(p!S^x_p)_{p\in\NN}$ and $\widehat{W}^x=(p!W^x_p)_{p\in\NN}$.

Let $\tau_1\in\hyperlink{omset}{\mathcal{W}}$ and $\tau_2\in\hyperlink{omset}{\mathcal{W}}$ be the weight functions coming from Theorem~\ref{theoEqualSpacesMatrixwithfactorialWeight} applied to $\sigma$ and $\omega$ respectively, so $\mathcal{A}_{\{\widehat{\Sigma}\}}(S_\gamma)= \mathcal{A}_{\{\tau_1\}}(S_\gamma)$ and $\mathcal{A}_{\{\widehat{\Omega}\}}(S_\gamma)= \mathcal{A}_{\{\tau_2\}}(S_\gamma)$.

Then, for every $l>0$ there exists $l_1>0$ such that there exists a linear and continuous map
\begin{equation}\label{coroExtensionOperatorsWeightsequ}
E^{\tau_1,\tau_2}_{l}:\;\;\;\lambda\in\Lambda_{\tau_1,l}\mapsto f_{\lambda}\in \mathcal{A}_{\tau_2,l_1}(S_{\gamma})
\end{equation}
such that for all $\lambda\in\Lambda_{\tau_1,l}$ one has $\mathcal{B}\circ E^{\tau_1,\tau_2}_l(\lambda)=\mathcal{B}(f_{\lambda})=\lambda$. Thus we have shown that $\mathcal{B}(\mathcal{A}_{\{\tau_2\}}(S_{\gamma}))\supseteq\Lambda_{\{\tau_1\}}$.
\end{corollary}

\demo{Proof}
Let $\mathcal{T}^i:=\{T^{i,x}: x>0\}$ be the weight matrix associated with the weight function $\tau_i$, i.e. $T^{i,x}_p:=\exp\left(\frac{1}{x}\varphi^{*}_{\tau_i}(xp)\right)$ for each $x>0$ and $p\in\NN$, $i=1,2$.

We may apply \eqref{equaEqualitySpacesWeightFunctionMatrix} in order to deduce that
$\mathcal{A}_{\{\widehat{\Sigma}\}}(S_\gamma)= \mathcal{A}_{\{\mathcal{T}^1\}}(S_\gamma)$ and $\mathcal{A}_{\{\widehat{\Omega}\}}(S_\gamma)= \mathcal{A}_{\{\mathcal{T}^2\}}(S_\gamma)$ and, as it has already been remarked in \cite[Corollary 7.6]{sectorialextensions}, we get $\mathcal{T}^1\hyperlink{Mroumapprox}{\{\approx\}}\widehat{\Sigma}$, $\mathcal{T}^2\hyperlink{Mroumapprox}{\{\approx\}}\widehat{\Omega}$.

For the rest of the proof we follow \cite[Corollary 7.6]{sectorialextensions} and use for the extension Theorem \ref{theoExtensionOperatorsMatrixmixed}.

\qed\enddemo

\begin{remark}\label{coroExtensionOperatorsWeightsrem}
Of course, by the assumption on the weights in Corollary \ref{coroExtensionOperatorsWeights}, we could apply \cite[Corollary 7.6]{sectorialextensions} directly to the weights $\omega$ and/or $\sigma$ because both have \hyperlink{om1}{$(\omega_1)$} and which is equivalent to have $\gamma(\omega),\gamma(\sigma)>0$ as shown in \cite[Corollary 2.14]{firstindexpaper}. Then \cite[Corollary 7.6]{sectorialextensions} states that for every $l>0$ there exists $l_1>0$ such that there exists a linear and continuous map
\begin{equation}\label{coroExtensionOperatorsWeightsremequ}
E^{\tau_i}_{l}:\;\;\;\lambda\in\Lambda_{\tau_i,l}\mapsto f_{\lambda}\in \mathcal{A}_{\tau_i,l_1}(S_{\gamma}),
\end{equation}
with $0<\gamma<\gamma(\sigma)$ for $i=1$ and $0<\gamma<\gamma(\omega)$ for $i=2$. Since in the mixed setting we are interested in the situation where $\sigma$ is (much) larger than $\omega$, \eqref{coroExtensionOperatorsWeightsremequ} yields in particular an extension $\Lambda_{\tau_1,l}\subseteq\Lambda_{\tau_2,l}\mapsto f_{\lambda}\in \mathcal{A}_{\tau_2,l_1}(S_{\gamma})$, $0<\gamma<\gamma(\omega)$. This should be compared with \eqref{coroExtensionOperatorsWeightsequ}, where we have such a mixed extension for all $0<\gamma<\gamma(\sigma,\omega)$ and since $\gamma(\sigma,\omega)\ge\gamma(\omega)$ (see Lemma \ref{indicescomparison}) even here we can have a situation which is not covered by \cite[Corollary 7.6]{sectorialextensions} since $\mathcal{A}_{\tau_2,l_1}(S_{\gamma'})\subseteq\mathcal{A}_{\tau_2,l_1}(S_{\gamma})$ for $\gamma\le\gamma'$.
\end{remark}

Finally we treat the mixed weight sequence situation and which is generalizing \cite[Remark 7.9]{sectorialextensions}. Note that condition \hyperlink{mg}{$(\on{mg})$} on the smaller weight sequence has also been used in \cite{ChaumatChollet94} and \cite[Theorem 5.10]{mixedramisurj}.

\begin{theorem}\label{theoExtensionOperatorssequencemixed}
Let $M,N\in\hyperlink{LCset}{\mathcal{LC}}$ be given such that $\mu_p\le C\nu_p$, $M$ satisfies \hyperlink{mg}{$(\on{mg})$} and finally $\gamma(M,N)>0$. We put $\widehat{M}=(p!M_p)_{p\in\NN}$ and $\widehat{N}=(p!N_p)_{p\in\NN}$ and let $0<\gamma<\gamma(M,N)(=\gamma(\omega_M,\omega_N))$.

Then, there exists a constant $k_1>0$ such that for every $h>0$, one can construct a linear and continuous map
$$
E^{M,N}_h:\;\;\lambda\in\Lambda_{\widehat{M},h}\mapsto f_{\lambda}\in \mathcal{A}_{\widehat{N},k_1h}(S_{\gamma})
$$
such that for every $\lambda\in\Lambda_{\widehat{M},h}$ one has $\mathcal{B}\circ E^{M,N}_h(\lambda)=\mathcal{B}(f_{\lambda})=\lambda$. Thus we have shown $\mathcal{B}(\mathcal{A}_{\{\widehat{N}\}}(S_{\gamma}))\supseteq\Lambda_{\{\widehat{M}\}}$.
\end{theorem}

\demo{Proof}
We apply Theorem \ref{theoExtensionOperatorsMatrixmixed} to $\sigma\equiv\omega_M$ and $\omega\equiv\omega_N$. Recall that by the assumptions on the weight sequences we have $0<\gamma(M,N)=\gamma(\omega_M,\omega_N)$, see Lemma \ref{mglemma1}. Moreover the matrix $\Sigma$ associated with $\omega_M$ is constant, see $(iii)$ in Lemma \ref{assofuncproper} and Remark \ref{importantremark}, and hence the matrix $\widehat{\Sigma}$ is constant, too. Theorem \ref{theoExtensionOperatorsMatrixmixed} provides now for all $h_1>0$ an extension map $E_{h_1}:\lambda\in\Lambda_{\widehat{S}^{1/8},h_1}\mapsto f_{\lambda}\in \mathcal{A}_{\widehat{W}^1,k_0h_1}(S_{\gamma})$, $k_0>0$ not depending on $h_1$.

We know that $S^1\equiv M$ (e.g. see the proof of \cite[Thm. 6.4]{testfunctioncharacterization}) and so $S^{1/8}\hyperlink{approx}{\approx} M\Leftrightarrow \widehat{S}^{1/8}\hyperlink{approx}{\approx}\widehat{M}$. Hence $\Lambda_{\widehat{S}^{1/8}}=\Lambda_{\widehat{M}}$, more precisely there exists some $D\ge 1$ such that for all $h>0$ we get $\Lambda_{\widehat{M},h}\subseteq\Lambda_{\widehat{S}^{1/8},Dh}$ and so we have shown the desired statement with universal $k_1=k_0D$.
\qed\enddemo

\section{Mixed extension results with only one fixed weight}\label{localweightsequence}

\subsection{Extension results where the weight sequence/function defining the function space is fixed}\label{localweightsequencesect1}
Using the properties of the index $\mu(N)$ and the construction of the ramified descendant of Section \ref{ordersofquasi} we can now prove the following variant of Theorem \ref{theoExtensionOperatorssequencemixed}.

\begin{theorem}\label{theoExtensionOperatorssequencemixedrem}
Let $N\in\hyperlink{LCset}{\mathcal{LC}}$ be given with $\mu(N)>0$ and let $0<r<\mu(N)$. Assume that \eqref{descendantmgequ} holds true for $N^{1/r}$. Then there does exist $L\in\hyperlink{LCset}{\mathcal{LC}}$ having \hyperlink{mg}{$(\on{mg})$} such that for each $0<\gamma<r$ we get: There exists a constant $k_1>0$ such that for every $h>0$, one can construct a linear and continuous map
\begin{equation}\label{theoExtensionOperatorssequencemixedremequ}
E^{L,N}_h:\;\;a\in\Lambda_{\widehat{L},h}\mapsto f_{a}\in\mathcal{A}_{\widehat{N},k_1h}(S_{\gamma}),
\end{equation}
denoting $\widehat{L}=(p!L_p)_{p\in\NN}$ and $\widehat{N}=(p!N_p)_{p\in\NN}$. Thus we have shown that $\mathcal{B}(\mathcal{A}_{\{\widehat{N}\}}(S_{\gamma}))\supseteq \Lambda_{\{\widehat{L}\}}$.

The sequence $L$ is maximal among those $M\in\hyperlink{LCset}{\mathcal{LC}}$ satisfying $\mu_k\le C\nu_k$ and \hyperlink{gammarmix}{$(M,N)_{\gamma_r}$}.
\end{theorem}

The important difference between Theorem \ref{theoExtensionOperatorssequencemixed} and this result is that, of course, $L$ is depending here on given $r$.

\demo{Proof}
Let $0<\gamma<r<\mu(N)$ be given according to the requirements above. Then we consider the sequence $L^{N,r}$ defined via the descendant $S^{N,r}$ in \eqref{sequenceL}, see Section \ref{ordersofquasi}.

As seen there we have that $\hyperlink{gammarmix}{(L^{N,r},N)_{\gamma_r}}$ holds true and which proves $\gamma(L^{N,r},N)\ge r>\gamma$. Moreover $\lambda^{N,r}_k\le C\nu_k$ and since $N^{1/r}$ has \eqref{descendantmgequ}, Lemma \ref{descendantmg} yields \hyperlink{mg}{$(\on{mg})$} for $S^{N,r}$ and so for $L^{N,r}$, too.

Thus we can apply Theorem \ref{theoExtensionOperatorssequencemixed} to $M\equiv L^{N,r}$ and $N$ and $\gamma$ unchanged to obtain: There exists a constant $k_1>0$ such that for every $h>0$, one can construct a linear and continuous map
$E^{L^{N,r},N}_h:\;\;a\in\Lambda_{\widehat{L}^{N,r},h}\mapsto f_{a}\in\mathcal{A}_{\widehat{N},k_1h}(S_{\gamma})$ with $\widehat{L}^{N,r}=(p!L^{N,r}_p)_{p\in\NN}$ and so \eqref{theoExtensionOperatorssequencemixedremequ} follows by taking $L\equiv L^{N,r}$.
\qed\enddemo

\begin{remark}
Let $N\in\hyperlink{LCset}{\mathcal{LC}}$ be given with $\mu(N)>0$. If $N$ has in addition \hyperlink{mg}{$(\on{mg})$}, then each $S^{N,r}$ and $L^{N,r}$, $0<r<\mu(N)$, share this property, see $(iv)$ in Remark \ref{descendant}.

If $N$ does not have \hyperlink{mg}{$(\on{mg})$}, then \eqref{descendantmgequ} for $N^{1/r}$ can only hold true for values $\gamma(N)\le r<\mu(N)$: for $0<r<\gamma(N)$ the right-hand side in \eqref{descendantmgequ} is bounded by definition, whereas $\sup_{k\in\NN}\frac{(\nu_{2k})^{1/r}}{(\nu_k)^{1/r}}=\infty$ since $N^{1/r}$ does also not have \hyperlink{mg}{$(\on{mg})$}.
\end{remark}


Using $\mu(\omega)$ we can prove Theorem \ref{theoExtensionOperatorssequencemixedrem} for the weight function setting, so we have the following variant of Theorem \ref{theoExtensionOperatorsMatrixmixed}.

\begin{theorem}\label{theoExtensionOperatorsMatrixmixedrem}
Let $\omega$ be a normalized weight function with \hyperlink{om3}{$(\omega_3)$} and $\mu(\omega)>0$. Then for all $0<r<\mu(\omega)$ there does exist a normalized weight function $\sigma$ satisfying \hyperlink{om3}{$(\omega_3)$} such that for each $0<\gamma<r$ we get: There exists a constant $k_0>0$ such that for every $x>0$ and every $h>0$, one can construct a linear and continuous map
$$
E^{\sigma,\omega}_h:\;\;\lambda\in\Lambda_{\widehat{S}^x,h}\mapsto f_{\lambda}\in \mathcal{A}_{\widehat{W}^{8x},k_0h}(S_{\gamma})
$$
such that for every $\lambda\in\Lambda_{\widehat{S}^x,h}$ one has $\mathcal{B}\circ E^{\sigma,\omega}_h(\lambda)=\mathcal{B}(f_{\lambda})=\lambda$. Thus we have shown that $\mathcal{B}(\mathcal{A}_{\{\widehat{\Omega}\}}(S_{\gamma})\supseteq \Lambda_{\{\widehat{\Sigma}\}}$ (by using for $\widehat{\Sigma}$ and $\widehat{\Omega}$ the same notation as in Theorem \ref{theoExtensionOperatorsMatrixmixed}).

The function $\sigma$ is chosen minimal among those normalized weight functions $\tau$ satisfying \hyperlink{om3}{$(\omega_3)$}, $\tau\hyperlink{ompreceq}{\preceq}\omega$ (i.e. $\omega(t)=O(\tau(t))$) and enjoying \hyperlink{gammarfctmix}{$(\tau,\omega)_{\gamma_r}$}.
\end{theorem}

\demo{Proof}
According to this value $r>0$ given, we consider the weight $\kappa^{1/r}_{\omega^r}$ (see \eqref{gammamixedvsomegamfctequ}) and so \hyperlink{gammarfctmix}{$(\kappa^{1/r}_{\omega^r},\omega)_{\gamma_r}$} is valid. This weight enjoys all properties like $\omega$ except normalization (by definition). But normalization can be achieved w.l.o.g. by switching to an equivalent weight (redefining $\kappa^{1/r}_{\omega^r}$ near $0$, e.g. see \cite[Remark 1.2 $(b)$]{BonetBraunMeiseTaylorWhitneyextension}) and which will be denoted by $\sigma$.

Thus $\gamma(\sigma,\omega)\ge r>\gamma$ and we can apply Theorem \ref{theoExtensionOperatorsMatrixmixed} to these weights $\sigma$ and $\omega$ and the value $\gamma$ and conclude.
\qed\enddemo

\begin{remark}
Due to \eqref{Qvsomegaequ} one could try to restate Theorem \ref{theoExtensionOperatorsMatrixmixedrem} by applying Theorem \ref{theoExtensionOperatorssequencemixedrem} to $N\equiv W^x$, $x>0$ arbitrary. However, once chosen $\gamma<\mu(\omega)=\mu(W^x)$ in Theorem \ref{theoExtensionOperatorsMatrixmixedrem} we obtain an extension for another weight function $\sigma$ such that moving the index $x$ we are staying in the same weight matrix associated with $\sigma$ by the precise choice $x\mapsto 8x$. So here we can take some uniform choice for all sequences in $\Omega$ (by obtaining a weight matrix not depending on given $x$) and which is not following by applying Theorem \ref{theoExtensionOperatorssequencemixedrem}.
\end{remark}

\begin{remark}\label{Javipredecessordestroyer}
Naturally one might ask what happens in the dual situation, that is, fixing the weight sequence or weight function that controls the derivatives at the origin. However, in this case the inverse construction concerning the descendant, called the {\itshape predecessor}, see \cite[Remark 4.3]{whitneyextensionweightmatrix}, does not provide any new information, since the bounds for the opening are the same as those that are known for the one level extension theorem.
\end{remark}

\appendix
\section{An example for having $\gamma(M,N),\gamma(\omega_M,\omega_N)>0$ but $\gamma(\omega_M)=\gamma(M)=\gamma(N)=\gamma(\omega_N)=0$}\label{examplesfromlangenbruch}

The aim is to construct an explicit example such that Theorem \ref{theoExtensionOperatorsMatrixmixed} resp. Theorem \ref{theoExtensionOperatorssequencemixed} are valid but no further known extension result construct can be applied and so the situation in this present work is not covered by the approach from \cite{sectorialextensions} (and also not from \cite{Thilliezdivision}, \cite{sectorialextensions1}). We will construct a pair of sufficient weight sequences $M$ and $N$ and apply Theorem \ref{theoExtensionOperatorssequencemixed} to $M$ and $N$ directly, whereas Theorem \ref{theoExtensionOperatorsMatrixmixed} is applied to their associated weights $\omega_M\equiv\sigma$ and $\omega_N\equiv\omega$.\vspace{6pt}

First recall that \hyperlink{gamma1}{$(\gamma_1)$} can be written in a different form. We say that $M$ has \hypertarget{beta1}{$(\beta_1)$}, if
$$\exists\;Q\in\NN_{>0}:\;\liminf_{p\rightarrow\infty}\frac{\mu_{Qp}}{\mu_p}>Q,$$
and \hypertarget{beta3}{$(\beta_3)$} (see \cite{dissertation} and \cite{BonetMeiseMelikhov07}), if
$$\exists\;Q\in\NN_{>0}:\;\liminf_{p\rightarrow\infty}\frac{\mu_{Qp}}{\mu_p}>1.$$
By \cite[Proposition 1.1]{petzsche} condition \hyperlink{beta1}{$(\beta_1)$} is equivalent to \hyperlink{gamma1}{$(\gamma_1)$} for log-convex $M$ and for this proof condition \hyperlink{mnq}{$(\text{nq})$}, which is a general assumption in \cite{petzsche}, was not necessary. We have that $M$ satisfies \hyperlink{beta1}{$(\beta_1)$} if and only if $m$ satisfies \hyperlink{beta3}{$(\beta_3)$}, for more precise information concerning this relation we refer to the recent work \cite{firstindexpaper}.\vspace{6pt}

The construction of the sequences is based on a generalization of \cite[Example 3.3]{Langenbruch89} and we introduce $M$ by using the quotient sequence $(\mu_k)_{k\ge 0}$ (with $\mu_0:=1$ and set $M_k:=\prod_{i=0}^k\mu_i$).

\begin{lemma}\label{generalizelangenbruch}
Let $\gamma>1$ be given, then there exists a sequence $M\in\hyperlink{LCset}{\mathcal{LC}}$ such that
\begin{itemize}
\item[$(i)$] $M$ does satisfy \hyperlink{mnq}{$(\on{nq})$}, more precisely $\mu_k\ge k^{\gamma}$ for all $k\in\NN$ and so even \hyperlink{mnqr}{$(\on{nq}_{\gamma-\varepsilon})$} holds true for any $\varepsilon>0$ (small),

\item[$(ii)$] $\mu_k\le k^{\gamma(2\gamma-1)}$ for all $k\in\NN$,

\item[$(iii)$] $M$ does not satisfy \hyperlink{beta3}{$(\beta_3)$} or equivalently $\widehat{M}=(p!M_p)_{p\in\NN}$ does not satisfy \hyperlink{beta1}{$(\beta_1)$} (and consequently $M$ is not strongly nonquasianalytic too),

\item[$(iv)$] $M$ does satisfy \hyperlink{mg}{$(\on{mg})$}.
\end{itemize}
\end{lemma}

\demo{Proof}
Let $\gamma>1$ be given and for convenience we introduce also numbers $\alpha$ and $\beta$ by
$$\gamma>1,\hspace{20pt}\beta:=2\gamma(>2),\hspace{20pt}\alpha:=\beta-1=2\gamma-1(>\gamma>1).$$
In \cite[Example 3.3]{Langenbruch89} the choices $\gamma=2$, $\beta=4$ and $\alpha=3$ have been considered (but there $(ii)$ has not been shown).

We define recursively two nondecreasing sequences $(c_n)_{n\ge 1}$ and $(d_n)_{n\ge 1}$ (of natural numbers) as follows and $\lfloor a\rfloor$ shall denote the lower integer part of a given real number $a>0$. We put
$$c_1:=1,\hspace{20pt}d_n:=\lfloor c_n^{\alpha/\gamma}\rfloor+1,\hspace{20pt}c_{n+1}:=\lfloor d_n^{\gamma}\rfloor+1,$$
and introduce the weight sequence as follows:
\begin{equation*}
\mu_k:=
\begin{cases}
c_n^{\alpha}\hspace{15pt}\text{for}\;c_n\le k\le\lfloor c_n^{\alpha/\gamma}\rfloor=d_n-1\\
\frac{k^{\beta}}{d_n^{\gamma}}\hspace{14pt}\text{for}\;d_n\le k\le\lfloor d_n^{\gamma}\rfloor=c_{n+1}-1.\\
\end{cases}
\end{equation*}

\qed\enddemo

$(i)$ and $(ii)$ together tell us that $M$ is lying between two Gevrey sequences. By proving $(i)$ one can verify that $\mu(M)=\liminf_{p\rightarrow\infty}\frac{\log(\mu_p)}{\log(p)}=\gamma$ holds true.\vspace{6pt}

A slight variation of Lemma \ref{generalizelangenbruch} yields the following.

\begin{lemma}\label{generalizelangenbruch0}
Let $\gamma>1$ be given, then there exists a sequence $M\in\hyperlink{LCset}{\mathcal{LC}}$ such that
\begin{itemize}
\item[$(i)$] $M$ does satisfy \hyperlink{mnq}{$(\on{nq})$}, more precisely $\mu_k\ge k^{\gamma}$ for all $k\in\NN$ and so even \hyperlink{mnqr}{$(\on{nq}_{\gamma-\varepsilon})$} holds true for any $\varepsilon>0$ (small),

\item[$(ii)$] $\mu_k\le k^{2\gamma^2}$ for all $k\in\NN$,

\item[$(iii)$] $M$ does not satisfy \hyperlink{beta3}{$(\beta_3)$} or equivalently $\widehat{M}=(p!M_p)_{p\in\NN}$ does not satisfy \hyperlink{beta1}{$(\beta_1)$} (and consequently $M$ is not strongly nonquasianalytic too),

\item[$(iv)$] $M$ does not satisfy \hyperlink{mg}{$(\on{mg})$}.
\end{itemize}
\end{lemma}

\demo{Proof}
We use the same definition for $(\mu_k)_k$ as before. Let $\gamma>1$ be given and set now
$$\gamma>1,\hspace{20pt}\beta:=2\gamma+1(>3),\hspace{20pt}\alpha:=\beta-1=2\gamma(>\gamma>1).$$
In fact any choice $\beta>2\gamma$ and $\alpha:=\beta-1$ would be working for the following proof.

\qed\enddemo

Again by construction $\mu(M)=\liminf_{k\rightarrow\infty}\frac{\log(\mu_k)}{\log(k)}=\gamma$ holds true.\vspace{6pt}

Using Lemma \ref{generalizelangenbruch} we can now underline the importance of Theorem \ref{theoExtensionOperatorsMatrixmixed} and in particular of Theorem \ref{theoExtensionOperatorssequencemixed} as follows.

\begin{theorem}\label{generalizelangenbruch1}
There do exist sequences $M$ and $N$ satisfying all requirements from Theorem \ref{theoExtensionOperatorssequencemixed} but such that $\gamma(M)=\gamma(N)=\gamma(\omega_M)=\gamma(\omega_N)=0$. Moreover we can achieve $\gamma(M,N)$ to be as large as desired.
\end{theorem}

\demo{Proof}
We define $M$ and $N$ in terms of their quotients $(\mu_p)_p$ and $(\nu_p)_p$ coming from Lemma \ref{generalizelangenbruch} with parameters $\gamma'$ and $\gamma$ respectively and we require that
\begin{equation}\label{gammarequirement}
\gamma'(2\gamma'-1)\le\gamma\Leftrightarrow\frac{\gamma'(2\gamma'-1)-\gamma}{r}\le 0.
\end{equation}
Choosing $1<\gamma'<\gamma$ subject to \eqref{gammarequirement} it is straightforward to prove \hyperlink{gammarmix}{$(M,N)_{\gamma_r}$} for all $0<r<\gamma$. This implies $\gamma(M,N)\ge\gamma>0$ and since $\mu(N)=\gamma$ we get equality by Lemma \ref{gammasmallerthanq}. Since $\gamma>1$ can be chosen arbitrary large, $\gamma(M,N)=\gamma$ can be as large as desired. Moreover one can easily verify $\mu_k\le\nu_k$, hence Lemma \ref{mglemma1} implies $\gamma(\omega_M,\omega_N)=\gamma(M,N)=\gamma$.\vspace{6pt}

But $\gamma(M)=\gamma(N)=0$ holds true: \hyperlink{beta1}{$(\beta_1)$} or equivalently \hyperlink{gamma1}{$(\gamma_1)$} is violated for both sequences $\widehat{M}=(p!M_p)_p$ and $\widehat{N}=(p!N_p)_p$ (by property $(iii)$), and so $\gamma(M)=\gamma(N)=0$. And this is equivalent to having $\gamma(\omega_{M})=\gamma(\omega_{N})=0$, because both sequences have \hyperlink{mg}{$(\on{mg})$}, for a proof see \cite[Section 4]{firstindexpaper}.
In particular we have seen that neither $\widehat{M}\in\hyperlink{SRset}{\mathcal{SR}}$ nor $\widehat{N}\in\hyperlink{SRset}{\mathcal{SR}}$ and by the characterizations shown in \cite{petzsche}, not any to $\widehat{M}$ or $\widehat{N}$ equivalent sequence $L$ can belong to class \hyperlink{SRset}{$\mathcal{SR}$}.

\qed\enddemo

Let $M$ and $N$ denote the sequences from Theorem \ref{generalizelangenbruch1} above with parameters $\gamma'>1$ and $\gamma>1$ subject to \eqref{gammarequirement}. Then by applying Theorem \ref{theoExtensionOperatorssequencemixed} for any given $0<\delta<\gamma(M,N)$ there is $k_1>0$ such that for every $h>0$ there exists a continuous linear extension map
$$E^{M,N}_h:\;\;\lambda\in\Lambda_{\widehat{M},h}\mapsto f_{\lambda}\in \mathcal{A}_{\widehat{N},k_1h}(S_{\delta}).$$
Thus we have shown $\mathcal{B}(\mathcal{A}_{\{\widehat{N}\}}(S_{\delta}))\supseteq \Lambda_{\{\widehat{M}\}}$.\vspace{6pt}

This kind of extension result is not covered by the theory developed by the authors in \cite{sectorialextensions}. More precisely \cite[Theorem 7.4]{sectorialextensions} fails since $\gamma(\omega_M)=\gamma(\omega_N)=0$ and also the mixed setting from \cite[Section 7.1]{sectorialextensions} cannot be applied, neither to $M$ nor to $N$ directly.\vspace{6pt}

Note that both $M$ and $N$ have \hyperlink{mg}{$(\on{mg})$}, thus both matrices associated with $\omega_M$ and $\omega_N$ are constant, see $(iii)$ in Lemma \ref{assofuncproper} and Remark \ref{importantremark}. Now let $M$ and $N$ be the sequences constructed in Lemma \ref{generalizelangenbruch0} with parameters $\gamma'$ and $\gamma$ respectively and here we require that
\begin{equation*}\label{gammarequirementnew}
2(\gamma')^2\le\gamma.
\end{equation*}
Again it is straightforward to check that \hyperlink{gammarmix}{$(M,N)_{\gamma_r}$} holds true for all $0<r<\gamma$ and which implies $\gamma(M,N)\ge\gamma>0$. Since $\mu(N)=\gamma$ we again have $\gamma(M,N)=\gamma$ and by having $\mu_p\le\nu_p$, Lemma \ref{mglemma1} yields $\gamma(\omega_M,\omega_N)\ge\gamma(M,N)=\gamma$.\vspace{6pt}

But here neither $M$ nor $N$ does satisfy \hyperlink{mg}{$(\on{mg})$} and we cannot apply Theorem \ref{theoExtensionOperatorssequencemixed} directly. But Theorem \ref{theoExtensionOperatorsMatrixmixed} applied to $\sigma\equiv\omega_M$ and $\omega\equiv\omega_N$ with $\Sigma$ denoting the matrix associated with $\omega_M$ and $\Omega$ the matrix associated with $\omega_N$, yields now the following extension result:\vspace{6pt}

For any given $0<\delta<\gamma(\omega_M,\omega_N)$ there exists a constant $k_0>0$ such that for every $x>0$ and every $h>0$, one can construct a linear and continuous extension map
$$
E^{\omega_M,\omega_N}_h:\;\;\lambda\in\Lambda_{\widehat{S}^x,h}\mapsto f_{\lambda}\in \mathcal{A}_{\widehat{W}^{8x},k_0h}(S_{\gamma}).
$$
Hence we have shown $\mathcal{B}(\mathcal{A}_{\{\widehat{\Omega}\}}(S_{\gamma}))\supseteq\Lambda_{\{\widehat{\Sigma}\}}$.\vspace{6pt}

{\itshape Note:} Here, since \hyperlink{mg}{$(\on{mg})$} is avoided, the associated matrices $\Sigma$ and $\Omega$, and hence $\widehat{\Sigma}$ and $\widehat{\Omega}$, are nonconstant but $S^1\equiv M$ and $W^1\equiv N$ (see the arguments given in the proof of Theorem \ref{theoExtensionOperatorssequencemixed}). Also in this situation we get $\gamma(M)=\gamma(N)=0$ (by property $(iii)$) but it is not clear if also $\gamma(\omega_M)=\gamma(\omega_N)=0$ (see \cite[Corollary 4.6]{firstindexpaper}). In such a situation, even if $\gamma(\omega_M)>0$ and/or $\gamma(\omega_N)>0$, we obtain new information, see Remark \ref{coroExtensionOperatorsWeightsrem} above.\vspace{6pt}

As mentioned in the introduction and in Remark \ref{omega1loss} we have that starting directly with a Braun-Meise-Taylor weight function $\omega$ with $\gamma(\omega)=0$ we do not have \hyperlink{om1}{$(\omega_1)$} (as shown in \cite[Corollary 2.14]{firstindexpaper}). Hence a basic assumption in the whole theory of ultradifferentiable functions defined in terms of $\omega$, is violated from the very beginning.\vspace{6pt}

\textbf{Acknowledgements}: The first two authors are partially supported by the Spanish Ministry of Economy, Industry and Competitiveness under the project MTM2016-77642-C2-1-P. The first author has been partially supported by the University of Valladolid through a Predoctoral Fellowship (2013 call) co-sponsored by the Banco de Santander. The third author is supported by FWF-Project J~3948-N35, as a part of which he has been an external researcher at the Universidad de Valladolid (Spain) for the period October 2016 - December 2018.\par

\bibliographystyle{plain}
\bibliography{Bibliography}

\vskip1cm

\textbf{Affiliation}:\\
J.~Jim\'{e}nez-Garrido, J.~Sanz:\\
Departamento de \'Algebra, An\'alisis Matem\'atico, Geometr{\'\i}a y Topolog{\'\i}a, Universidad de Valladolid\\
Facultad de Ciencias, Paseo de Bel\'en 7, 47011 Valladolid, Spain.\\
Instituto de Investigaci\'on en Matem\'aticas IMUVA\\
E-mails: jjjimenez@am.uva.es (J.~Jim\'{e}nez-Garrido), jsanzg@am.uva.es (J. Sanz).
\\
\vskip.5cm
G.~Schindl:\\
Fakult\"at f\"ur Mathematik, Universit\"at Wien,
Oskar-Morgenstern-Platz~1, A-1090 Wien, Austria.\\
E-mail: gerhard.schindl@univie.ac.at.


\end{document}